\documentclass[12pt]{article}
\usepackage{amssymb}
\usepackage{latexsym}
\usepackage{amsmath}
\usepackage{amscd}
\usepackage{color }
\newfont{\cyr}{wncyr10}
\newfont{\cyb}{wncyb10}

\def\sA {{\mathcal A}}
\def\sB {{\mathcal B}}

\def\sD {{\mathcal D}}
\def\sE {{\mathcal E}}

\def\sH {{\mathcal H}}
\def\sI {{\mathcal I}}
\def\sJ {{\mathcal J}}
\def\sK {{\mathcal K}}
\def\sL {{\mathcal L}}
\def\sM {{\mathcal M}}

\def\sP {{\mathcal P}}

\def\sS {{\mathcal S}}
\def\sT {{\mathcal T}}
\def\sU {{\mathcal U}}
\def\sV {{\mathcal V}}
\def\sW {{\mathcal W}}

\def\od {\mathrm{od}}

\def\int {\mathrm{int}}
\def\cl {\mathrm{cl}}

\def\ps {\mathrm{ps}}

\def\nwd {\mathrm{nwd}}

\newtheorem{thm}{Theorem}[section]

\newtheorem{theorem}[thm]{Theorem}
\newtheorem{corollary}[thm]{Corollary}
\newtheorem{lemma}[thm]{Lemma}

\newtheorem{definition}[thm]{Definition}

\newtheorem{discussion}[thm]{Discussion}

\newtheorem{notation}[thm]{Notation}
\newtheorem{remark}[thm]{Remark}
\newtheorem{remarks}[thm]{Remarks}
\newtheorem{shan}[thm]{Standing Hypotheses and Notation}

\newtheorem{question}[thm]{Question}

\newtheorem{notation-definition}[thm]{Notation and Definitions}

\begin{document}


\title{Tychonoff Expansions with Prescribed Resolvability Properties}

\author{W.W. Comfort\footnote{Email: wcomfort@wesleyan.edu; 
Address: Department of Mathematics and Computer Science,
           Wesleyan University, Wesleyan Station,
           Middletown, CT 06459; phone: 860-685-2632; FAX 860-685-2571}, 
Wanjun Hu\footnote{Corresponding author, Email: Wanjun.Hu@asurams.edu; 
Address: Department of Mathematics and Computer Science, Albany State University,
            Albany, GA 31705; phone: 229-886-4751}}

\maketitle

\begin{abstract}
The recent literature offers examples, specific and
hand-crafted, of Tychonoff spaces (in ZFC) which respond negatively to
these questions, due respectively to Ceder and Pearson
(1967) and to Comfort and Garc\'{\i}a-Ferreira (2001): (1)~Is every
$\omega$-resolvable space maximally resolvable?
(2)~Is every maximally resolvable space
extraresolvable? Now using the method of ${\mathcal{KID}}$ expansion, the
authors show that {\it every} suitably restricted
Tychonoff topological space $(X,\sT)$ admits a larger Tychonoff
topology (that is, an``expansion")
witnessing such failure. Specifically the authors show in ZFC that if
$(X,\sT)$ is a maximally resolvable Tychonoff space
with $S(X,\sT)\leq\Delta(X,\sT)=\kappa$, then
$(X,\sT)$ 
has Tychonoff expansions $\sU=\sU_i$ ($1\leq i\leq5$), with
$\Delta(X,\sU_i)=\Delta(X,\sT)$ and $S(X,\sU_i)\leq\Delta(X,\sU_i)$,
such that $(X,\sU_i)$ is:
($i=1$)~$\omega$-resolvable but not maximally resolvable;
($i=2$)~[if $\kappa'$ is regular, with
$S(X,\sT)\leq\kappa'\leq\kappa$] $\tau$-resolvable for all
$\tau<\kappa'$, but not
$\kappa'$-resolvable; 
($i=3$)~maximally resolvable, but not extraresolvable;
($i=4$)~extraresolvable, but not maximally resolvable;
($i=5$)~maximally resolvable and extraresolvable, but not
strongly extraresolvable.


{\sl MSC Primary 05A18, 03E05, 54A10; Secondary 03E35, 54A25, 05D05}

\end{abstract}


\section{Introduction, Definitions and Notation}

Our principal interest is in Tychonoff spaces, i.e., in completely regular,
Hausdorff spaces, and
all spaces $(X,\sT)$ hypothesized here, 
also all expansions
(refinements) of $\sT$ constructed, will
be Tychonoff topologies. The topological properties we consider, however, are
intelligible (a wonderful word in this context, borrowed from
Hewitt~\cite{hewa}) for arbitrary spaces, so in \ref{defprops} below, which
defines the properties we consider,
we impose no separation hypotheses.

\begin{notation}
{\rm
For $X$ a set and $\tau$ a cardinal, we set $[X]^\tau:=\{A\subseteq
X:|A|=\tau\}$. The symbols $[X]^{<\tau}$ and $[X]^{\leq\tau}$ are
defined analogously. 

The symbol $D(\tau)$ denotes the discrete space of cardinality $\tau$.

When $X=(X,\sT)$ is a space and $Y\subseteq X$, we denote by $(Y,\sT)$
the set $Y$ with the subspace topology inherited from $X$.

The symbols $w$ and $d$ denote {\it weight} and {\it density character},
respectively. For a space $X=(X,\sT)$, the {\it dispersion character}
$\Delta(X)$ is the smallest cardinal of an nonempty open subset of $X$,
and $\nwd(X)$, the {\it nowhere density number of} $X$, is

$\nwd(X):=\min\{|A|:A\subseteq X, \int_X\,\cl_X\,A\neq\emptyset\}$.\\
\noindent Evidently $\nwd(X)$
coincides with the {\it open density number of} $X$~\cite{comfhu03}
defined by

$\od(X):=\min\{d(U):\emptyset\neq U\in\sT\}$,\\
\noindent which has also been denoted $d_0(X)$~\cite{juhss2}.

As in \cite{comfhu07} and \cite{juhss}, a subset $D$ of a space $X=(X,\sT)$ is
$\tau$-{\it dense} in $X$ if $|D\cap U|\geq\tau$ whenever $\emptyset\neq
U\in\sT$. It is obvious that if $D$ is dense in a space $X$ with
$\nwd(X)\geq\tau$, then $D$ is $\tau$-dense in $X$.

$(X,\sT)$ is {\it crowded} if no
point of $X$ is isolated in the topology $\sT$. (This term,
introduced by van~Douwen~\cite{vd93b}, has been adopted subsequently by
many authors \cite{eck97}, \cite{hu}, \cite{juhss}. Others have called
such a space {\it dense-in-itself} \cite{comfhu04}.)

A family of nonempty  pairwise disjoint open subsets of $X=(X,\sT)$ is a
{\it cellular family}, and $S(X)$, the {\it Souslin number} of $X$, is

$S(X):=\min\{\kappa:$ no cellular $\sU\subseteq\sT$ satisfies $|\sU|=\kappa\}$.
}
\end{notation}

\begin{definition}\label{defprops}
{\rm
Let $X=(X,\sT)$ be a space. Then $X$ is

\begin{itemize}

\item[(i)] {\it resolvable} (Hewitt~\cite{hewa}) if it has two complementary
dense subsets;

\item[(ii)] {\it $\kappa$-resolvable} (Ceder~\cite{ceder64}) if there is a
family of $\kappa$-many pairwise disjoint dense subsets of $X$;

\item[(iii)] {\it maximally resolvable} (Ceder~\cite{ceder64}) if it is
$\Delta(X)$-resolvable;

\item[(iv)] {\it extraresolvable} (Malykhin~\cite{maly98}) if there is a family
$\sD$ of dense subsets, with $|\sD|\geq(\Delta(X))^+$, such that every
two elements of $\sD$ have interesection which is nowhere dense in $X$;
and

\item[(v)] {\it strongly extraresolvable}
(Comfort and Garc\'{\i}a-Ferreira~\cite{comfgarc98}, \cite{comfgarc01})
if there is a family
$\sD$ of dense subsets, with $|\sD|\geq(\Delta(X))^+$, such that
distinct $D_0$, $D_1\in\sD$ satisfy $|D_0\cap D_1|<\nwd(X)$.
\end{itemize}
}
\end{definition}

\begin{remark}
{\rm
In early versions of this manuscript, circulated privately to selected
colleagues, we were able to establish item ($i=4$) of the Abstract, even
its special case Theorem~\ref{i=4}, only under the additional assumption
that there exists a cardinal $\tau$ such that $\tau<\kappa<2^\tau$.
Indeed, although we had shown in \cite{comfhu07} the existence of
extraresolvable Tychonoff spaces which are not maximally resolvable when
GCH fails, it was an unsolved problem whether such spaces exist
in ZFC. That question has been settled affirmatively by Juh\'{a}sz,
Shelah and Soukup~\cite{juh-sh-sou}. We are grateful to those authors
for furnishing us with a pre-publication copy of their work.
}
\end{remark}

\begin{definition}\label{DIndep}
{\rm
Let $\kappa\geq\omega$.

(a) A partition $\sB$ of $\kappa$ is a $\kappa$-partition if each
$B\in\sB$ satisfies $|B|=\kappa$;

(b) a family $\sB=\{\sB_t:t\in T\}$ of partitions
$\sB=\{\sB^\alpha_t:\alpha<\kappa_t\}$
of $\kappa$ is
$\tau$-{\it independent} (with $1\leq\tau\leq\kappa$)
if $|\bigcap_{t\in F}\,B^{f(t)}_t|\geq\tau$ for each
$F\in[T]^{<\omega}$ and $f\in\Pi_{t\in F}\,\kappa_t$.

(c) a family $\sB=\{\sB_t:t\in T\}$ of
indexed partitions $\sB_t=\{B_t^\alpha:\alpha<\kappa_t\}$ (with
$2\leq\kappa_t\leq\kappa$ for each $t\in T$)
{\it separates points} [resp., {\it separates small sets}] if for distinct
$x,x'\in\kappa$ there are
$\sB_t\in\sB$ and (distinct) $\alpha,\alpha'<\kappa_t$ such that $x\in
B_t^\alpha$ and $x'\in B_t^{\alpha'}$ [resp., for disjoint
$S,S'\in[\kappa]^{<\kappa}$ there are
$\sB_t\in\sB$ and (distinct) $\alpha,\alpha'<\kappa_t$ such that
$S\subseteq B_t^\alpha$ and $S'\subseteq B_t^{\alpha'}$].
}
\end{definition}

It is obvious that any partition in a $\kappa$-independent family (of
partitions of $\kappa$) is necessarily a $\kappa$-partition.

\begin{discussion}\label{disc1}
{\rm
Given a point-separating family ${\sB}$ as in
Definition~\ref{DIndep}, we denote by $\sT_{\sB}$ the smallest
topology on $\kappa$ in which
each set $B_t^\alpha\in\sB_t\in\sB$ is open; clearly each such
$B_t^\alpha$ is $\sT_\sB$-closed, and
$\{\bigcap_{t\in F}\,B_t^{f(t)}:F\in[T]^{<\omega},f\in\Pi_{t\in
F}\,\kappa_t\}$ is a basis for $\sT_\sB$.
(This is a Hausdorff topology since $\sB$ separates points of $\kappa$,
hence is a Tychonoff topology since it has a clopen basis.)
The evaluation map
$e_\sB:(\kappa,\sT_\sB)\rightarrow\Pi_{t\in T}\,D(\kappa_t)$ given by
$$
(e_\sB x)_t=\alpha\textit{ if }x\in B_t^\alpha~~~~~(x\in\kappa,
t\in T, \alpha<\kappa_t)
$$
\noindent is a homeomorphism from $(\kappa,\sT_\sB)$ onto a
subspace $X$ of the Tychonoff space $K:=\Pi_{t\in T}\,D(\kappa_t)$. That
$X:=e_\sB[\kappa]$ is dense in $K$
follows from the fact that $\sB$ is $1$-independent. Conversely,
given $K=\Pi_{t\in T}\,D(\kappa_t)$ with $|T|=2^\kappa$ and
with $2\leq\kappa_t\leq\kappa$ for each $t\in T$, the
Hewitt-Marczewski-Pondiczery theorem (cf.~\cite{engel}(2.3.15),
\cite{comfneg74}(\S3 and Notes)) gives
a dense set $X\subseteq K$ such that $|X|=\kappa$, and then the family 
${\sB}:=\{\sB_t:t\in T\}$ with
$\sB_t:=\{\pi_t^{-1}(\{\alpha\}\cap X:\alpha<\kappa_t\}$
is a $1$-independent family of partitions of $\kappa$
(the set $\kappa$ here being identified with the subspace $X$ of $K$).
One may ensure that each $\sB_t\in\sB$ is a $\kappa$-partition by the following
device (here we argue much as in \cite{comfhu04}(3.8) and
\cite{comfhu07}(1.5)): Give each space $D(\kappa_t)$ the structure
of a topological
group, so that $K$ is a topological group, let $X^*$ be dense in $K$
with $|X^*|=\kappa$, and with $\langle X^*\rangle$ the subgroup of $K$
generated by $X^*$ let $X$ be the union of $\kappa$-many cosets of
$\langle X^*\rangle$ in $K$. Then
$B^\alpha_t:=\pi_t^{-1}(\{\alpha\})\cap X$
satisfies $|B_t^\alpha|=\kappa$ for each $\alpha<\kappa_t$, $t\in T$;
indeed more generally each basic open set $U$ in $X$ (of the form
$U=(\bigcap_{i=1}^n\,\pi_{t_i}^{-1}(\{\alpha_i\}))\cap X$, with
$\alpha_i<\kappa_{t_i}$, $n<\omega$) satisfies $|U|=\kappa$, so
the family $\sB$ is even $\kappa$-independent, and $\Delta(X)=\kappa$.

The correspondence ${\sB}\leftrightarrow X$
just described is of Galois type in the sense that 
when dense $X\subseteq K=\Pi_{t\in T}\,D(\kappa_t)$ is given with
$|X|=\kappa$ and the family ${\sB}=\{\sB_t:t\in T\}$ is defined, then
$e_{\sB}:(\kappa,\sT_{\sB})\rightarrow K$ satisfies $e_{\sB}[\kappa]=X$.
}
\end{discussion}

In this paper in this context, $T$ and $\{\kappa_t:t\in T\}$ being given, we
use the notations $(\kappa,\sT_\sB)$, $(X,\sT_{\sB})$ and $e_\sB[\kappa]$
interchangeably.

The point-separating
family described in Discussion~\ref{disc1}
may be chosen to separate small sets in a strong
sense. Lemma~\ref{sss}, which exploits a trick
introduced by Eckertson~\cite{eck97}
in a related context, strengthens a statement given in our
works \cite{comfhu03} and
\cite{comfhu04}(3.3(b)).
When reference is made, in Lemma~\ref{sss} and later, to a partition
$\{T(\lambda):\lambda\in\Lambda\}$ of $T$, the trivial (one-cell)
partition is not excluded.

\begin{lemma}\label{sss}
Let $\kappa\geq\omega$ and $|T|=2^\kappa$, and for $t\in T$ let
$2\leq\kappa_t\leq\kappa$. Let $\{T(\lambda):\lambda\in\Lambda\}$ be a
partition of $T$, with each $|T(\lambda)|=2^\kappa$.
Then there is a $\kappa$-independent
family $\sI=\{\sI_t:t\in T\}$ of partitions of $\kappa$,
with $|\sI_t|=\kappa_t$ for each
$t\in T$,
such that for every ordered pair
$(S,S')$ of disjoint elements of $[\kappa]^{<\kappa}$
and for every $\lambda\in\Lambda$ there are
infinitely many $t\in T(\lambda)$ such that $S\subseteq I^0_t$ and
$S'\subseteq I^1_t$.
\end{lemma}

Proof. Let ${\sB}=\{\sB_t:t\in T\}$ be
a point-separating $\kappa$-independent family of partitions of
$\kappa$ with $|T|=2^\kappa$ and with $|\sB_t|=\kappa_t$
for each $t\in T$, 
as given in Discussion~\ref{disc1}.
For $\lambda\in\Lambda$ let $\{T(\lambda,\xi):\xi<2^\kappa\}$ be a
partition of $T(\lambda)$ with each $|T(\lambda,\xi)|=\omega$,
and using $|[\kappa]^{<\kappa}|\leq2^\kappa$ let
$\{(S_\xi,S'_\xi):\xi<2^\kappa\}$ list all ordered pairs of disjoint
members of $[\kappa]^{<\kappa}$ (with repetitions permitted). Then
define $\sI=\{\sI_t:t\in T\}$ with
$\sI_t=\{I_t^\alpha:\alpha<\kappa_t\}$
as follows: 
if $t\in T(\lambda,\xi)$, then

$I_t^0=(B_t^0\cup S_\xi)\backslash S'_\xi$,
$I_t^1=(B_t^1\cup S'_\xi)\backslash S_\xi$, and
$I_t^\alpha=B_t^\alpha\backslash(S_\xi\cup S'_\xi)$
for $2\leq\alpha<\kappa_t$.\\
\noindent Then each $\sI_t$ is a partition of $\kappa$, and
since

$B^\alpha_t\triangle I^\alpha_t\in[\kappa]^{<\kappa}$\hfill (*)\\
\noindent for each $t\in T$ and $\alpha<\kappa_t$ with $\sB_t$ a
$\kappa$-partition, so also is each $\sI_t$ a $\kappa$-partition.
Further for each pair $(S,S')=(S_\xi,S'_\xi)$ we have $S\subseteq I^0_t$
and $S'\subseteq I^1_t$ for each
$t\in T(\lambda,\xi)\in[T(\lambda)]^\omega$, as required.~$\square$

\begin{definition}
{\rm
With $\{\kappa_t:t\in T\}$ and $\{T(\lambda):\lambda\in\Lambda\}$ given
as in Lemma~\ref{sss},
a $\kappa$-independent family $\sI$ of partitions of
$\kappa$ with the additional property given there is a {\it
strong small-set-separating} family of partitions which {\it respects} the
partition
$\{T(\lambda):\lambda\in\Lambda\}$ of $T$.
}
\end{definition}

\begin{remark}
{\rm
Clearly a $\kappa$-independent family $\{\sI_t:t\in T\}$
of partitions of $\kappa$, if it respects some partition
$\{T(\lambda):\lambda\in\Lambda\}$ of $T$, also respects the trivial
(one-cell) partition. Usually in this paper
we apply Lemma~\ref{sss} only in the
context of the trivial partition; in what follows, if no explicit
reference is made to the partition which a strong small-set-separating
family of $\kappa$-partitions respects, we intend by default the
trivial partition.
}
\end{remark}

The following theorem augments, simplifies
and extends arguments given in our works \cite{comfhu04}(3.8) and
\cite{comfhu07}(1.6). As usual when a point-separating family $\sI$ of
partitions of $\kappa$ is given, we do not distinguish notationally
between $\kappa$ and the space $X:=e_\sI[\kappa]\subseteq K=\Pi_{t\in
T}\,D(\kappa_t)$, nor between a set $I^\alpha_t\in\sI_t\in\sI$ and its
image $e_\sI[I^\alpha_t]$ in $X$.

\begin{theorem}\label{nwdX=kappa}
Let $\kappa\geq\omega$ and $|T|=2^\kappa$, and for $t\in T$ let
$2\leq\kappa_t\leq\kappa$. Then there is a $\kappa$-independent
family $\sI=\{\sI_t:t\in T\}$ of partitions of $\kappa$
with the strong small-set-separating property, and
with $|\sI_t|=\kappa_t$ for each
$t\in T$,
such that the space

$X:=e_\sI[\kappa]\subseteq K:=\Pi_{t\in T}\,D(\kappa_t)$\\
\noindent has these properties:

{\rm (a)} $X$ is dense in $K$;

{\rm (b)} $X$ is $\kappa$-resolvable;

{\rm (c)} $|X|=\Delta(X)=\nwd(X)=\kappa$; and

{\rm (d)} each $S\in[X]^{<\kappa}$ is closed and discrete in $X$.
\end{theorem}

Proof. Let $\overline{T}:=T\cup\{\overline{t}\}$ with $\overline{t}\notin
T$, and set $\kappa_{\overline{t}}:=\kappa$.
Apply Lemma~\ref{sss} with $\{\overline{T}\}$ the one-cell
partition of $\overline{T}$:
There is a $\kappa$-independent
family $\overline{\sI}=\{\sI_t:t\in \overline{T}\}$
of $\kappa$-partitions of $\kappa$ with
the strong small-set-separating property,
with $|\sI_t|=\kappa_t$ for each $t\in\overline{T}$ (in particular, with
$|\sT_{\overline{t}}|=\kappa_{\overline{t}}=\kappa$). 
By the argument given in Discussion~\ref{disc1} the set
$X:=e_{\sI}[\kappa]$ is dense in $K:=\Pi_{t\in T}\,D(\kappa_t)$, so (a)
is proved. For $F\in[T]^{<\omega}$ and $f\in\Pi_{t\in F}\,\kappa_t$ and
each $I^\alpha_{\overline{t}}$ (with
$\alpha<\kappa_{\overline{t}}=\kappa$) we have

$|(\bigcap_{t\in F}\,I^{f(t)}_t)\cap
I^\alpha_{\overline{t}}|=\kappa$\hfill (*)\\
\noindent since the family
$\overline{\sI}$ is $\kappa$-independent. Relation (*) shows that each set
$e_\sI[I^\alpha_{\overline{t}}]$ is dense in $X$ (thus proving (b)), and it
shows also that $|X|=\Delta(X)=\kappa$.

Since $X$ is a crowded space, every closed, discrete subspace of $X$ is
nowhere dense; so the relation $\nwd(X)=\kappa$ will follow from (d).
Given $S\in[\kappa]^{<\kappa}$ and
$x\in \kappa\backslash S$, there is $t\in T$ such that
$x\in I^0_t$ and $S\subseteq I^1_t$;
since $I^0_t$ and $I^1_t$ are disjoint and clopen in $X$, we conclude
that $S$ is closed. Similarly if $x\in S\in[\kappa]^{<\kappa}$
there is $t\in T$ such that
$x\in I^0_t$ and $S\backslash\{x\}\subseteq I^1_t$, so $I^0_t\cap
S=\{x\}$; it follows that $S$ is discrete.~$\square$

\begin{remarks}
{\rm

(a) In earlier work \cite{comfhu07}
by a different argument we have demonstrated the
existence of a $\kappa$-resolvable dense subset $X$ of some spaces of the form
$\Pi_{t\in T}\,D(\kappa_t)$ with $|T|=2^\kappa$, even with
$|X|=\Delta(X)=\nwd(X)=\kappa$. (See also
\cite{comfhu03}(5.3 and 5.4) for similar results.)
The argument of Theorem~\ref{nwdX=kappa}
is preferable, both because of its simplicity and because
it gives in concrete form a family $\sI$ for which
$X=e_{\sI}[\kappa]$; this latter feature is essential in
the proof of Lemma~\ref{S_er} below.

(b) The case in Definition~\ref{DIndep} in which
there is $\lambda\in[2,\kappa]$ such that $\kappa_t=\lambda$ for all
$t\in T$, together with passage in that case from
${\sB}$
to the space $(\kappa,\sT_{\sB})=(X,\sT_{\sB})$, has been used by
many authors in connection with resolvability questions
\cite{vd93b}, \cite{comfhu03},
\cite{comfhu04}, \cite{juhss}, \cite{comfhu07}.
} 
\end{remarks}

\section{The ${\mathcal{KID}}$ Expansion: Transfer from $\sT$ to
$\sT_{\mathcal{KID}}$}

Here we explain and develop further the techniques originating in
\cite{huthesis}, \cite{hu}. In broad terms the goal, given a
crowded Tychonoff
space $(X,\sT)$, is to augment (``expand") the topology $\sT$ to a larger
crowded Tychonoff
topology $\sT_{\mathcal{KID}}$ in such a way that certain specified
$\sT$-dense subsets of $X$
remain $\sT_{\mathcal{KID}}$-dense,
while certain other subsets of $X$ become closed and discrete in the
topology $\sT_{\mathcal{KID}}$.

In Definition~\ref{defKID}, the transition from $\sT$ to the
$\sT_{\mathcal{KID}}$-open sets $W^\alpha_t$ is effected {\it via} the
intermediate sets $H^\alpha_t$. Their definition depends on the
hypothesized dense array $\sD$ and the $\kappa$-independent family
$\sI$, but not on the family $\sK$.

The following notation is as in \cite{comfhu04}(3.2).

\begin{notation}\label{X(S)}
{\rm
Let $X$ be a set with $|X|=\kappa\geq\omega$, and let
$\sD=\{D^\gamma_\eta:\gamma<\tau,\eta<\kappa\}$ be a partition of $X$
with $1\leq\tau\leq\kappa$. Then for $S\subseteq\kappa$ the set
$X(S)\subseteq X$ is defined by

$X(S):=\bigcup\{D^\gamma_\eta:\gamma<\tau,\eta\in S\}$.
}
\end{notation}

\begin{definition}\label{defKID}
{\rm
Let $(X,\sT)$ be a crowded Tychonoff space with $|X|=\kappa\geq\omega$,
fix nonempty $Z\subseteq X$,
and let $\sI=\{\sI_t:t\in Z\times2^\kappa\}$ be a
point-separating $\kappa$-independent family of
partitions of $\kappa$ with $\sI_t=\{I^\alpha_t:\alpha<\kappa_t\}$,
$2\leq\kappa_t\leq\kappa$ for each $t\in Z\times2^\kappa$.
Let $1\leq\tau\leq\kappa$ and
$\sD=\{D^\gamma_\eta:\gamma<\tau,\eta<\kappa\}$ be a
partition of $X$, and
for $t\in Z\times2^\kappa$ and $\alpha<\kappa_t$ set

$H^\alpha_t:=X(I^\alpha_t)=\bigcup\{D^\gamma_\eta:\gamma<\tau,\eta\in
I^\alpha_t\}$.

Let $\sK=\{K_\xi:\xi<2^\kappa\}\subseteq\sP(Z)$,
and for $t=(x,\xi)\in Z\times2^\kappa$ and $\alpha<\kappa_t$ define
$W^\alpha_t$ as follows:

If $K_\xi=\emptyset$, then $W^\alpha_t = H^\alpha_t$.

If $K_\xi\neq\emptyset$, then

$W^0_t=(H^0_t\cup K_\xi)\backslash\{x\}$,

$W^1_t=(H^1_t\backslash K_\xi\big)\cup\{x\}$, and

$W^\alpha_t=H^\alpha_t\backslash(K_\xi\cup\{x\})$ for
$2\leq\alpha<\kappa_t$.

For each $t\in Z\times 2^\kappa$ set

$\sH_t:=\{H^\alpha_t:\alpha<\kappa_t\}$
and $\sW_t:=\{W^\alpha_t:\alpha<\kappa_t\}$,\\
\noindent and set

$\sH:=\{\sH_t:t\in Z\times2^\kappa\}$, and
$\sW:=\{\sW_t:t\in Z\times2^\kappa\}$.

\noindent Then

$\sT^{\mathcal {ID}}$
is the smallest topology on $X$ such that
$\sT\subseteq\sT^{\mathcal{ID}}$ and each
$\sH_t\subseteq\sT^{\mathcal{ID}}$, and

$\sT_{\mathcal{KID}}$, the ${\mathcal{KID}}$ {\it expansion} of $\sT$, is
the smallest topology on $X$ such that
$\sT\subseteq\sT_{\mathcal{KID}}$ and each
$\sW_t\subseteq\sT_{\mathcal{KID}}$.
}
\end{definition}

\begin{remarks}\label{KIDremarks}
{\rm
(a) The indexings $\sD=\{D^\gamma_\eta:\gamma<\tau,\eta<\kappa\}$
and $\sI=\{\sI_t:t\in Z\times2^\kappa\}$ in Definition~\ref{defKID} are
faithful. No such restriction is imposed on the indexing
$\sK=\{K_\xi:\xi<2^\kappa\}$. Indeed in some of the applications we
will have $K_\xi=\emptyset$ for many $\xi<2^\kappa$.

(b) For $t\in Z\times2^\kappa$
the family $\sH_t$ is a partition of
$X$ into $\sT^{\mathcal{ID}}$-open subsets,
so each $H^\alpha_t$ is $\sT^{\mathcal{ID}}$-clopen. 
Similarly, since for $t\in Z\times2^\kappa$ the family
$\sW_t$ is a partition
of $X$ into $\sT_{\mathcal{KID}}$-open sets, also each
$W^\alpha_t$ is $\sT_{\mathcal{KID}}$-clopen. 
It then follows, as is required of every topology hypothesized or
constructed in this paper, that:

(c) Each space of the form $(X,\sT^{\mathcal{ID}})$, and
each space of the form $(X,\sT_{\mathcal{KID}})$, is a Tychonoff space.

(d) the topology $\sT_{\mathcal{KID}}$ depends not only on the families
$\sK$, $\sI$, and $\sD$, but also on the choice of the
nonempty set $Z\subseteq X$. Our notation does not reflect that fact. No
confusion with ensue, indeed in (nearly) all the applications we take
$Z=X$. Briefly in Theorem~\ref{usingS_er} we will invoke the
general theory in the special case $|Z|=1$.

}
\end{remarks}

To avoid irrelevancies we gave Definition~\ref{defKID} in
uncluttered language, but in fact we will use the expansion
$\sT_{\mathcal{KID}}$ only when the following additional
conditions are satisfied. Except when noted otherwise,
we assume these henceforth throughout
this Section. Furthermore when families $\sI$, $\sD$ and $\sK$ have been
constructed or hypothesized and $I^\alpha_t\in\sI_t\in\sI$, it is
understood that the sets $H^\alpha_t$ and $W^\alpha_t$ are defined as
in Definition~\ref{defKID}.

\begin{shan}\label{shan1}~~

{\rm

(1)~$|X|=\Delta(X,\sT)=\kappa$;

(2)~the indexed family
$\sD$ is a dense partition of
$(X,\sT)$, and $D^\gamma:=\bigcup_{\eta<\kappa}\,D^\gamma_\eta$
for $\gamma<\tau$;

(3)~the family $\sI=\{\sI_t:t\in
Z\times2^\kappa\}$ has the strong small-set-separating property;

(4)~if $F\in[2^\kappa]^{<\omega}$ then
$\bigcup_{\xi\in F}\,K_\xi\in\sK$; and

(5) $\xi<2^\kappa$,
$\gamma<\tau\Rightarrow\int_{(D^\gamma,\sT^{\mathcal{ID}})}(K_\xi\cap
D^\gamma)=\emptyset$.
}
\end{shan}

\begin{lemma}\label{KIDprops}
{\rm [With the conventions of {\rm \ref{defKID}} and {\rm \ref{shan1}}.]}\\
Fix $\gamma<\tau$ and $\xi<2^\kappa$. Then

{\rm (a)} $K_\xi$ is closed in
$(Z,\sT_{\mathcal{KID}})$; 

{\rm (b)} $(K_\xi,\sT_{\mathcal{KID}})$ is discrete; and

{\rm (c)} if $\emptyset\neq U\in\sT$, $H=\bigcap_{t\in F}\,H_t^{f(t)}$
and $W=\bigcap_{t\in F}\,W_t^{f(t)}$
with $F\in[Z\times2^\kappa]^{<\omega}$ and $f\in\Pi_{t\in F}\,\kappa_t$,
then $|D^\gamma\cap U\cap H|=|D^\gamma\cap U\cap W|=\kappa$.
\end{lemma}

Proof. 
(a) If $x\in Z\backslash K_\xi$ then with $t:=(x,\xi)$ we have $x\in
W^1_t\in\sT_{\mathcal{KID}}$ and $W^1_t\cap K_\xi=\emptyset$.

(b) If $x\in K_\xi$ then with $t:=(x,\xi)$ we have $W^1_t\in\sT_{\mathcal
KID}$ and $W^1_t\cap K_\xi=\{x\}$.

(c) Let $I:=\bigcap_{t\in F}\,I_t^{f(t)}$. Since $\sI$ is
$\kappa$-independent we have $|I|=\kappa$. For each $\eta\in I$ the set
$D^\gamma\cap H$ contains the set $D^\gamma_\eta$; since the sets
$D^\gamma_\eta$ ($\eta\in I$) are
pairwise disjoint, each dense in $(X,\sT)$, we have

$\kappa=|X|\geq|D^\gamma\cap U\cap H|\geq|I|=\kappa$.\hfill(*)\\
It remains to show that $|D^\gamma\cap U\cap W|=\kappa$. First, set

$K:=\bigcup_{(x,\xi)\in F}\,K_\xi$ and
$L:=\bigcup_{(x,\xi)\in F}\,(K_\xi\cup\{x\})$,\\
\noindent and note from (4) and (5) of \ref{shan1} that
$\int_{(D^\gamma,\sT^{\mathcal{ID}})}(D^\gamma\cap K)=\emptyset$, hence
also

$\int_{(D^\gamma,\sT^{\mathcal{ID}})}(D^\gamma\cap
L)=\emptyset$\hfill(**)\\
\noindent (since $(D^\gamma,\sT^{\mathcal{ID}})$ is crowded).

Now let $A:=(D^\gamma\backslash L)\cap(U\cap H)$. Since
$D^\gamma\cap U\cap W\supseteq A$, it suffices to show $|A|=\kappa$. If
$A\in[X]^{<\kappa}$ then, writing
$S:=\{\eta<\kappa:A\cap D^\gamma_\eta\neq\emptyset\}$,
we have $|S|\leq|A|<\kappa$, so by \ref{shan1}(3) there is
$\widetilde{t}\in(Z\times2^\kappa)\backslash F$ such that $S\subseteq
I^0_{\widetilde{t}}$; then $S\cap I^1_{\widetilde{t}}=\emptyset$
and hence $A\cap H^1_{\widetilde{t}}=\emptyset$. Then with

$\widetilde{f}:=f\cup\{(\widetilde{t},1)\}\in\Pi_{t\in
F\cup\{\widetilde{t}\}}\,\kappa_t$ and

$\widetilde{H}:=\bigcap_{t\in F\cup\{\widetilde{t}\}}\,H_t^{f(t)}=H\cap
H^1_{\widetilde{t}}\in\sH$\\
\noindent we have
$\emptyset=A\cap H^1_{\widetilde{t}}=
(D^\gamma\backslash L)\cap(U\cap H)\cap H^1_{\widetilde{t}}=
(D^\gamma\backslash L)\cap(U\cap\widetilde{H})$ and hence\\
$D^\gamma\cap L\supseteq(D^\gamma\cap L)\cap(U\cap\widetilde{H})
=\emptyset\cup[(D^\gamma\cap L)\cap(U\cap\widetilde{H})]$\\
$~~~~~~~~~=[(D^\gamma\backslash L)\cap(U\cap\widetilde{H})]
\cup[(D^\gamma\cap L)\cap(U\cap\widetilde{H})]$\\
$~~~~~~~~~=D^\gamma\cap U\cap\widetilde{H}$.\hfill(***)\\
\noindent By (*) applied with $\widetilde{H}$ replacing $H$, the set
$D^\gamma\cap U\cap\widetilde{H}$ is a nonempty
$\sT^{\mathcal{ID}}$-open subset of $D^\gamma$,
so (***) contradicts (**).~$\square$

\begin{corollary}\label{delta=delta}
{\rm [With the conventions of {\rm \ref{defKID}} and {\rm \ref{shan1}}.]}

{\rm (a)} $(D^\gamma,\sT^{\mathcal{ID}})$ is crowded, and $D^\gamma$
is dense in
$(X,\sT^{\mathcal{ID}})$;

{\rm (b)} $(D^\gamma,\sT_{\mathcal{KID}})$ is crowded, and $D^\gamma$
is dense in
$(X,\sT_{\mathcal{KID}})$; and

{\rm (c)} $\Delta(X,\sT^{\mathcal{ID}})=\Delta(X,\sT_{\mathcal
KID})=\Delta(X,\sT)=\kappa$.
\end{corollary}

Proof. The inequalities
$\Delta(X,\sT^{\mathcal{ID}})\leq\Delta(X,\sT)=\kappa$ and
$\Delta(X,\sT_{\mathcal{KID}})\leq\Delta(X,$ $\sT)$ $=\kappa$ of (c) follow
from the inclusions $\sT\subseteq\sT^{\mathcal{ID}}$ and
$\sT\subseteq\sT_{\mathcal{KID}}$, and all else is immediate
from Lemma~\ref{KIDprops}.~$\square$

It is easily seen that each infinite (Hausdorff) space $(X,\sT)$
contains an infinite cellular family, hence satisfies
$S(X,\sT)\geq\omega^+$. According to a result of Erd\H{o}s and
Tarski~\cite{erdtarski43} (see also \cite{comfneg74}(3.5),
\cite{comfneg82}(2.14)) every
infinite Souslin number is regular. That allows us to compute
exactly numbers of the form $S(X,\sT_{\mathcal{KID}})$ in terms of
the number $S(X,\sT)$ and the
family $\{\kappa_t:t\in Z\times2^\kappa\}$.

\begin{lemma}\label{S(KID)}
{\rm [With the conventions of {\rm \ref{defKID}} and {\rm \ref{shan1}}.]}\\
$S(X,\sT_{\mathcal{KID}})$ is the
smallest regular cardinal $\kappa'$ such that

{\rm (i)} $\kappa'\geq S(X,\sT)$, and

{\rm (ii)} $t\in Z\times2^\kappa\Rightarrow\kappa'\geq\kappa_t^+$.
\end{lemma}

Proof. From $\sT\subseteq\sT_{\mathcal{KID}}$ follows $S(X,\sT)\leq
S(X,\sT_{\mathcal{KID}})$. Further for $t\in Z\times2^\kappa$ the family
$\{W^\alpha_t:\alpha<\kappa_t\}$ is cellular in $(X,\sT_{\mathcal{KID}})$,
so $S(X,\sT_{\mathcal{KID}})\geq\kappa_t^+$. Since $S(X,\sT_{\mathcal
KID})$ is regular by the cited theorem of Erd\H{o}s and Tarski, we have
$S(X,\sT_{\mathcal{KID}})\geq\kappa'$.

Suppose now that $\{U_\zeta\cap W_\zeta:\zeta<\kappa'\}$ is a faithfully
indexed cellular family of $\sT_{\mathcal{KID}}$-basic open subsets of
$X$; here $U_\zeta\in\sT$ and $W_\zeta=\bigcap_{t\in
F_\zeta}\,W_t^{f_\zeta(t)}$ with
$F_\zeta\in[Z\times2^\kappa]^{<\omega}$, $f_\zeta\in\Pi_{t\in
F_\zeta}\,\kappa_t$, $W_t^{f_\zeta(t)}\in\sW$. Since
$\{F_\zeta:\zeta<\kappa'\}$
is a family of finite sets indexed (not necessarily faithfully) by the
regular cardinal $\kappa'$, there are $A\in[\kappa']^{\kappa'}$ and a
set $F$ such that $F_{\zeta_0}\cap F_{\zeta_1}=F$ for every pair
$\{\zeta_0,\zeta_1\}\in[A]^2$. (See \cite{comfneg74} or \cite{comfneg82}
or \cite{juhasz80} for proofs and bibliographic commentary on this
theorem, its special cases and generalizations.) Since $|F|<\omega$ and
$f_\zeta(t)<\kappa_t<\kappa'$ for each $\zeta\in A$ and $t\in F$, there is
$B\in[A]^{\kappa'}$ such that $f_{\zeta_0}(t)=f_{\zeta_1}(t)$ for all
$\zeta_0,\zeta_1\in B$ and $t\in F$. We define

$f:F_{\zeta_0}\cup F_{\zeta_1}\rightarrow\bigcup_{t\in F_{\zeta_0}\cup
F_{\zeta_1}}\,\kappa_t$\\
\noindent by

\[ f(t)= \left\{
\begin{array}
{r@{\quad \mathrm{if}\quad}l}
f_{\zeta_0}(t)=f_{\zeta_1}(t)       &   t\in F  \\
f_{\zeta_0}(t)   &  t\in F_{\zeta_0}\backslash F\\
f_{\zeta_1}(t)   &  t\in F_{\zeta_1}\backslash F\\
\end{array} \right\}.  \]
\noindent

\noindent (More succintly: $f=f_{\zeta_0}|F_{\zeta_0}\cup
f_{\zeta_1}|F_{\zeta_1}$.) Then since $S(X,\sT)\leq\kappa'=|B|$ there are
distinct $\zeta_0,\zeta_1$ (henceforth fixed) in $B$ such that
$U_{\zeta_0}\cap U_{\zeta_1}\neq\emptyset$. 

Then $H_{\zeta_0}\cap H_{\zeta_1}=\bigcap_{t\in F_0\cup F_1}\,H_t^{f(t)}$,
and (using (c) in Lemma~\ref{KIDprops}) we have

$\emptyset\neq(H_{\zeta_0}\cap H_{\zeta_1})\cap(U_{\zeta_0}\cap
U_{\zeta_1})\in\sT^{\mathcal{ID}}$.\\
\noindent Now choose and fix $\gamma<\tau$, and (arguing much as in the
proof of Lemma~\ref{KIDprops}(c)) set

$K:=\bigcup_{(x,\xi)\in F_{\zeta_0}\cup F_{\zeta_1}}\,K_\xi$ and
$L:=\bigcup_{(x,\xi)\in F_{\zeta_0}\cup F_{\zeta_1}}\,(K_\xi\cup\{x\})$;\\
\noindent then $K\in\sK$ by \ref{shan1}(4) and $D^\gamma\backslash K$ is
dense in the crowded space ($D^\gamma,\sT^{\mathcal{ID}})$ by \ref{shan1}(5),
so $D^\gamma\backslash L$ is also dense in
$(D^\gamma,\sT^{\mathcal{ID}})$,
hence also in $(X,\sT^{\mathcal{ID}})$ by
Corollary~\ref{delta=delta}(a). Then

$(D^\gamma\backslash L)\cap(H_{\zeta_0}\cap H_{\zeta_1})\cap(U_{\zeta_0}\cap
U_{\zeta_1})\neq\emptyset$,

\noindent so

$(D^\gamma\backslash L)\cap(W_{\zeta_0}\cap W_{\zeta_1})\cap(U_{\zeta_0}\cap
U_{\zeta_1})\neq\emptyset$,\\
\noindent contrary to the condition
$(W_{\zeta_0}\cap U_{\zeta_0})\cap(W_{\zeta_1}\cap
U_{\zeta_1})=\emptyset$.~$\square$

\begin{discussion}
{\rm
The method of ${\mathcal{KID}}$ expansion was introduced in
\cite{huthesis} and was
used in \cite{hu} to give the existence, assuming Lusin's Hypothesis,
of $\omega$-resolvable Tychonoff spaces which are not maximally resolvable.
The present authors have used the method subsequently \cite{comfhu04},
\cite{comfhu07} to find and construct explicit spaces with some of the
properties given in the Abstract. Arguments with some
similar features were found independently and exploited by Juh\'asz,
Szentmiklossy, and Soukup~\cite{juhss}; we note that
\cite{juhss} was submitted to the journal of record before
\cite{comfhu07} was submitted, furthermore the date of publication of
\cite{juhss} precedes that of \cite{comfhu07}.

The principal thrust of the present paper is this: Not only do specific
spaces (constructed as in \cite{huthesis}, \cite{hu}, \cite{comfhu04},
\cite{juhss}, \cite{comfhu07}) exist with the properties
listed, but indeed every crowded Tychonoff space subject to minimal
necessary conditions admits such Tychonoff expansions.
}
\end{discussion}

\begin{definition}\label{defnew}
{\rm [With the conventions of {\rm \ref{defKID}},
but with $\sK$ not yet defined.]
Let $\sM=\{M_\xi:\xi<2^\kappa\}\subseteq\sP(Z)$ with $M_0=\emptyset$. Then
$\widetilde{\sM}=\{\widetilde{M_\xi}:\xi<2^\kappa\}$ is defined as
follows.

$\widetilde{M_0}=\emptyset$, and

\noindent if $0<\xi<2^\kappa$ and
$\widetilde{M_{\eta}}$ has been defined for all $\eta<\xi$
then

$\widetilde{M_{\xi}}=M_{\xi}$ if each set of the form

$(M_\xi\cup\widetilde{M_{\eta_0}}\cup\widetilde{M_{\eta_1}}\cup\cdots\cup
\widetilde{M_{\eta_m}})\cap D^\gamma$ ($m<\omega$, $\eta_i<\xi$, $\gamma<\tau$)

\noindent has empty interior in the space $(D^\gamma,\sT^{\mathcal{ID}})$,

$\widetilde{M_\xi}=\emptyset$ otherwise.
}
\end{definition}

\begin{lemma}\label{discrete}
Let $Y$ be a crowded (Hausdorff)
space and let $E=\bigcup_{i\leq m}\,E_i\subseteq Y$
with each $E_i$ discrete, $m<\omega$. Then $\int_Y\,E=\emptyset$.
\end{lemma}

Proof. This is clear when $m=0$. Suppose it holds for $m=k$ and let
$E=\bigcup_{i\leq k+1}\,E_i\subseteq Y$ with each $E_i$ discrete.
Suppose for a contradiction that there is $p\in\int_Y\,E$, say with 
$p\in E_{k+1}$, and find open $U\subseteq Y$ such that $U\cap
E_{k+1}=\{p\}$. Then $(U\cap\int_Y\,E)\cap E_{k+1}=\{p\}$, so
$\bigcup_{i\leq k}\,E_i$ contains the nonempty open set
$(U\cap\int_Y\,E)\backslash\{p\}$.~$\square$

\begin{theorem}\label{heredirr}
{\rm [With the conventions of {\rm \ref{defKID}} and
{\rm \ref{shan1}(1), (2), (3)}.]}\\
Let $\sM=\{M_\xi:\xi<2^\kappa\}=\sP(Z)$ and let
$\sK:=\widetilde{\sM}=\{\widetilde{M_\xi}:\xi<2^\kappa\}$.
Then

{\rm (a)} $\sK$ satisfies conditions {\rm (4)} and {\rm (5)} of
{\rm \ref{shan1}};

{\rm (b)} if $\overline{\xi}<2^\kappa$ and
$\int_{(D^\gamma,\sT_{\mathcal{KID}})}\,(M_{\overline \xi}\cap
D^\gamma)=\emptyset$ for all $\gamma<\tau$, then
$M_{\overline \xi}=\widetilde{M_{\overline \xi}}\in\sK$; and

{\rm (c)} each space $(D^\gamma\cap Z,\sT_{\mathcal{KID}})$ is
hereditarily irresolvable.
\end{theorem}

Proof. (a) is obvious.

(b) 
Fix ${\overline \xi}<2^\kappa$ and $\gamma<\tau$, and let
$\eta_0,\eta_1,\cdots\eta_m<\overline{\xi}$, $\emptyset\neq U\in\sT$ and
$H=\bigcap_{t\in F}\,H^{f(t)}_t$ with $F\in[Z\times2^\kappa]^{<\omega}$,
$f\in\Pi_{t\in F}\,\kappa_t$.
We must show that if
$\int_{(D^\gamma,\sT_{\mathcal{KID}})}(M_{\overline{\xi}}\cap
D^\gamma)=\emptyset$ for all $\gamma<\tau$, then

$\int_{(D^\gamma,\sT^{\mathcal{ID}})}((M_{\overline{\xi}}\cup\widetilde{M_{\eta_0}}\cup\widetilde{M_{\eta_1}}\cup\ldots\cup
\widetilde{M_{\eta_m}})\cap D^\gamma)=\emptyset$.\hfill(*)\\
\noindent Writing $W=\bigcap_{t\in F}\,W^{f(t)}_t$,
we have, since $D^\gamma\backslash M_{\overline{\xi}}$ is dense in
$(D^\gamma,\sT_{\mathcal{KID}})$ and
$\emptyset\neq U\cap W\in\sT_{\mathcal{KID}}$, that

$Y:=(D^\gamma\backslash M_{\overline \xi})\cap(U\cap W)$ is dense in
$((D^\gamma\cap(U\cap W)),\sT_{\mathcal{KID}})$.\\
\noindent Further
since $(D^\gamma\cap(U\cap W),\sT_{\mathcal{KID}})$ is crowded, its
dense subset $(Y,\sT_{\mathcal{KID}})$ is crowded. 

We have $W\backslash H\subseteq L:=\bigcup_{(x,\xi)\in
F}\,(K_\xi\cup\{x\})$, with $L$ the union of finitely many discrete
subsets of $(Z,\sT_{\mathcal{KID}})\subseteq(X,\sT_{\mathcal{KID}})$.
Each $\widetilde{M_{\eta_i}}\in\sK$ is also discrete in
$(Z,\sT_{\mathcal{KID}})\subseteq(X,\sT_{\mathcal{KID}})$, so from
Lemma~\ref{discrete} it follows that the set
$Y\backslash(\bigcup_{i\leq m}\,\widetilde{M_{\eta_i}}\cup L)$
remains dense in $(D^\gamma\cap U\cap W,\sT_{\mathcal{KID}})$, and (*)
follows.

(c) Suppose for some $\gamma_0<\tau$ there are $\xi_0<2^\kappa$ and
nonempty $S\subseteq D^{\gamma_0}\cap Z$ such that $M_{\xi_0}\subseteq S$ and
both $M_{\xi_0}$ and $S\backslash M_{\xi_0}$ are dense in
$(S,\sT_{\mathcal{KID}})$.
From $\int_{(S,\sT_{\mathcal{KID}})}M_{\xi_0}=\emptyset$
it follows that $\int_{(D^{\gamma_0},\sT_{\mathcal
KID})}M_{\xi_0}=\emptyset$, so
$\int_{(D^\gamma,\sT_{\mathcal
KID})}(M_{\xi_0}\cap D^\gamma)=\emptyset$ for each $\gamma<\tau$. From
(b) we then have $M_{\xi_0}=\widetilde{M_{\xi_0}}\in\sK$, so by
Lemma~\ref{KIDprops}(a) the set $M_{\xi_0}$ is closed in
$(Z ,\sT_{\mathcal{KID}})$ (hence in $(S,\sT_{\mathcal{KID}})$); this
contradicts the density in $(S,\sT_{\mathcal{KID}})$ of both $M_{\xi_0}$
and $S\backslash M_{\xi_0}$.~$\square$

\section{The ${\mathcal{KID}}$ Expansion: Applications}

We begin this Section by proving (the case $|X|=\Delta(X)$ of)
our principal theorem (cf. item
($i=1$) of the Abstract). The result is in the tradition of the several
papers listed in the Bibliography which respond to the Ceder-Pearson
question (Is there an $\omega$-resolvable space which is not maximally
resolvable?), but this has a different flavor: Not only can examples of
such spaces be constructed by {\it ad hoc} techniques, but indeed {\it
every} (suitably restricted) $\omega$-resolvable Tychonoff space admits
a Tychonoff expansion $\sU$ such that $(X,\sU)$ remains
$\omega$-resolvable but is not maximally resolvable. For remarks
intended to justify or to explain the special hypothesis
``$S(X,\sT)\leq|X|$" in Theorem~\ref{maintheorem}, see
Remark~\ref{aboutquestions} below, where it is noted that in some
settings where $S(X,\sT)\leq|X|$ fails, $\omega$-resolvability implies
maximal resolvability.

\begin{theorem}\label{maintheorem}
Let $X=(X, \sT)$ be a crowded,
$\omega$-resolvable Tychonoff space 
with $S(X,\sT)\leq|X|= \Delta(X,\sT)=\kappa$.
Then there is a
Tychonoff refinement $\sU$ of $\sT$ such that

{\rm (a)} $S(X,\sU)=S(X,\sT)$ and $\Delta(X,\sU)=\Delta(X,\sT)$;

{\rm(b)}  $(X, \sU)$ is 
$\omega$-resolvable;

{\rm (c)} $(X,\sU)$ is not maximally resolvable; and

{\rm (d)} $(X,\sU)$ is not $S(X,\sT)$-resolvable, if $(X,\sT)$ is
maximally resolvable.
\end{theorem}

Proof. If $(X,\sT)$ is not maximally resolvable the conditions are
satisfied with $\sU:=\sT$, so we
assume in what follows that $(X,\sT)$ is maximally resolvable.

Let $\sD=\{D^n_\eta:\eta<\kappa,n<\omega\}$ be a faithfully indexed dense
partition of $(X,\sT)$, and set $D^n:=\bigcup_{\eta<\kappa}\,D^n_\eta$
for $n<\omega$.
Take $Z=X$ in Definition~\ref{defKID} and
let $\sI=\{\sI_t:t\in X\times 2^\kappa\}$ be a
$\kappa$-independent family of partitions
$\sI_t$ of $X$ with the strong small-set-separating property given by
Lemma~\ref{sss}; for simplicity we take $\kappa_t=2=\{0,1\}$ for each $t\in
X\times2^\kappa$.

Let $\sM=\{M_\xi:\xi<2^\kappa\}=\sP(X)$, and define
$\sK:=\widetilde{\sM}$ as in Definition~\ref{defnew}.
We will show that
$\sU:=\sT_{\mathcal{KID}}$ is as required.

(a) The equality $\Delta(X,\sT_{\mathcal{KID}})=\Delta(X,\sT)$ is given by
Corollary~\ref{delta=delta},
while $S(X,$ $\sT_{\mathcal{KID}})$ $=S(X,\sT)$ is immediate from
Lemma~\ref{S(KID)} (using the regularity of $S(X,\sT)$ and the fact that
$\kappa_t<\omega<\omega^+\leq S(X,\sT)$ for each $t\in Z\times2^\kappa$).

(b) According to Corollary~\ref{delta=delta}(b), the disjoint sets $D^n$
($n<\omega$) are dense in $(X,\sT_{\mathcal{KID}})$.

(c) and (d)  Suppose there is a
family $\sE$ of pairwise disjoint dense subsets of
$(X,\sT_{\mathcal{KID}})$ such that $|\sE|=S(X,\sT)$.
Note then that for some $E\in\sE$ we have

$\int_{(D^n,\sT_{\mathcal{KID}})}(D^n\cap E)=\emptyset$ for each
$n<\omega$.\hfill(*)\\
\noindent (Indeed otherwise we may argue as in \cite{hu}(2.3),
\cite{comfhu04}, \cite{comfhu07}(3.1(c)): choosing for each $E\in\sE$ some
$n(E)<\omega$ such that

$\int_{(D^{n(E)},\sT_{\mathcal{KID}})}(D^{n(E)}\cap E)\neq\emptyset$,\\
\noindent we have from Lemma~\ref{S(KID)} and
the regularity of $S(X,\sT)=S(X,\sT_{\mathcal{KID}})$ that some (fixed)
$n<\omega$ satisfies

$\int_{(D^n,\sT_{\mathcal{KID}})}(D^n\cap E)\neq\emptyset$ for
$S(X,\sT_{\mathcal{KID}})$-many $E\in\sE$;\\
\noindent that gives $S(D^n,\sT_{\mathcal{KID}})>S(X,\sT_{\mathcal
KID})$,
which is impossible since $D^n$ is dense in $(X,\sT_{\mathcal{KID}})$.)

Then choosing $E\in\sE$ as in (*), we have from
Theorem~\ref{heredirr}(b) that $E\in\sK$, so $E$
is closed and discrete in the crowded space $(X,\sT_{\mathcal
KID})$ by Lemma~\ref{KIDprops}((a) and (b)). This contradicts
the density of $E$ in
$(X,\sT_{\mathcal{KID}})$.~$\square$

\begin{remark}
{\rm
The choice $\kappa_t<\kappa$ for all $t\in X\times2^\kappa$ in (the
proof of) Theorem~\ref{maintheorem} is essential. If $\kappa_t=\kappa$
is permitted for some $t$ then the refinement $\sU=\sT_{\mathcal{KID}}$
satisfies conditions (b) and (c), but as noted in the first paragraph of
the proof of Lemma~\ref{S(KID)} we would now have $S(X,\sT_{\mathcal
KID})=\kappa^+>S(X,\sT)$.
}
\end{remark}

As is indicated in its proof, Theorem~\ref{maintheorem} is of interest
only when the given space $(X,\sT)$ is maximally resolvable. So viewed,
the case $\kappa'=S(X,\sT)$ of the
following result (cf. item ($i=2$) of our Abstract) strengthens
and improves Theorem~\ref{maintheorem}.

\begin{theorem}\label{notkappa'resolv}
Let $X=(X, \sT)$ be a crowded, maximally resolvable Tychonoff space
and let $\kappa'$ be a regular cardinal such that
$S(X,\sT)\leq\kappa'\leq|X|= \Delta(X,\sT)=\kappa$.
Then there is a
Tychonoff refinement $\sU$ of $\sT$ such that

{\rm (a)} $S(X,\sU)=\kappa'$ and $\Delta(X,\sU)=\Delta(X,\sT)=\kappa$;

{\rm(b)}  $(X, \sU)$ is 
$\tau$-resolvable for each $\tau<\kappa'$; and

{\rm (c)} $(X,\sU)$ is not $\kappa'$-resolvable.
\end{theorem}

Proof. [Being $\kappa$-resolvable, the space $(X,\sT)$ is surely
$\kappa'$-resolvable, so in this case the topology $\sU$ will of
necessity be a strict refinement of $\sT$.]

Let $\Lambda$ be the set of all cardinals $\tau$ such that
$2\leq\tau<\kappa'$, and let
$\{\kappa_t:t\in T=X\times2^\kappa\}$ list the elements of $\Lambda$
with each $\tau\in\Lambda$ appearing $2^\kappa$-many times.
For $\tau\in\Lambda$ set $T(\tau):=\{t\in T:\kappa_t=\tau\}$. According
to Lemma~\ref{sss}, there is a strong small-set-separating family
$\kappa$-independent family $\sI=\{\sI_t:t\in X\times2^\kappa\}$ of
partitions of $\kappa$ which respects the partition
$\{T(\tau):\tau\in\Lambda\}$ of $T$.

We note that $\kappa'=\sup_{t\in T}\,\kappa_t^+$.

Let $\sD=\{D^n_\eta:n<\omega,\eta<\kappa\}$ be a dense partition
of $(X,\sT)$, and as usual set $D^n:=\bigcup_{\eta<\kappa}\,D^n_\eta$.

Take $\sK$ as in Theorem~\ref{heredirr} and set
$\sU:=\sT_{\mathcal{KID}}$ (with $Z=X$). We show that $\sU$ is as required.

(a) The equalities $\Delta(X,\sT_{\mathcal{KID}})=\Delta(X,\sT)$ and
$S(X,\sT_{\mathcal{KID}})=\kappa'$ are given by Corollary~\ref{delta=delta}(c)
and Lemma~\ref{S(KID)}, respectively.

(c) The argument showing that the space
$(X,\sT_{\mathcal{KID}})$ of Theorem~\ref{maintheorem}(c)
is not $S(X,\sT_{\mathcal
KID})$-resolvable (i.e., is not $\kappa'$-resolvable)
applies here {\it verbatim} to prove (c).

(b) Let $\sA=\{A_n:n<\omega\}$ be an arbitrary countable dense partition
of the space $(X,\sT_{\mathcal{KID}})$.
Fix $\tau<\kappa'$, let $t(n)$ ($n<\omega$) be a faithfully indexed
sequence from $X\times2^\kappa$ such that $\kappa_{t(n)}=\tau$ for each
$n<\omega$, and for $n<\omega$ and $\alpha<\tau$ set

$E^\alpha_n:=W^\alpha_{t(n)}\backslash\bigcup_{k<n}\,W^\alpha_{t(k)}$.\\
\noindent Each set $E^\alpha_n$ is nonempty, and by
Remark~\ref{KIDremarks}(b) each is $\sT_{\mathcal{KID}}$-clopen. Now define

$E^\alpha:=\bigcup_{n<\omega}(E^\alpha_n\cap A_n)$  $(\alpha<\tau$);\\
\noindent we will show that $\{E^\alpha:\alpha<\tau\}$ is a dense
partition of $(X,\sT_{\mathcal{KID}})$.

Suppose there is $x\in E^\alpha\cap E^{\alpha'}$ with
$\alpha,\alpha'<\tau$. Then there are $n,n'<\omega$ such that

$x\in (E^\alpha_n\cap A_n)\cap(E^{\alpha'}_{n'}\cap A_{n'})\subseteq
A_n\cap A_{n'}$,\\
\noindent so $n=n'$ and from
$x\in E^\alpha_n\cap E^{\alpha'}_n\subseteq W^\alpha_{t(n)}\cap
W^{\alpha'}_{t(n)}$ we have
$\alpha=\alpha'$, as required.

To see for (fixed) $\alpha<\tau$ that $E^\alpha$ is dense in
$(X,\sT_{\mathcal{KID}})$, let $U\cap W\in\sT_{\mathcal
KID}$ with $\emptyset\neq U\in\sT$ and with $W=\bigcap_{t\in F}\,W^{f(t)}_t$ with
$F\in[X\times2^\kappa]^{<\omega}$, $f\in\Pi_{t\in F}\,\kappa_t$. We
assume without loss of generality, replacing $W$ by a smaller set if
necessary, that some $t(n)\in F$; and further with $m:=\max\{n:t(n)\in
F\}$ that $n<m\Rightarrow t(n)\in F$. It suffices to show that
$(U\cap W)\cap E^\alpha_n\neq\emptyset$ for some $n$, for then (from the
density of $A_n$ in $(X,\sT_{\mathcal{KID}})$ and the fact that
$E^\alpha_n$ is open in $(X,\sT_{\mathcal{KID}})$) it will follow that

$(U\cap W)\cap E^\alpha\supseteq(U\cap W)\cap(E^\alpha_n\cap A_n)
=(U\cap W\cap E^\alpha_n)\cap A_n\neq\emptyset$.

\underline{Case 1}. Some $n\leq m$ satisfies $f(t(n))=\alpha$. Then,
choosing minimal such $n$, we have $\emptyset\neq U\cap W
\subseteq W\subseteq E^\alpha_n$, so
$(U\cap W)\cap E^\alpha_n=U\cap W\neq\emptyset$.

\underline{Case 2}. Case 1 fails. Then defining
$\widetilde{f}:=f\cup\{(t(m+1),\alpha)\}$
we have $W\cap W^\alpha_{t(m+1)}\subseteq E^\alpha_{t(m+1)}$,
and Lemma~\ref{KIDprops}(c) gives

$\emptyset\neq U\cap(W\cap W^\alpha_{t(m+1)})\cap E^\alpha_{t(m+1)}
\subseteq (U\cap W)\cap(E^\alpha_{m+1})$.~$\square$

\begin{remarks}\label{aboutkappa'-res}
{\rm
(a) According to Theorem~\ref{KIDprops} the family $\{D^n:n<\omega\}$ is
a dense partition of $(X,\sT_{\mathcal{KID}})$. We note that the
construction just given parlays an arbitrary countable dense partition
$\sA=\{A_n:n<\omega\}$ of $(X,\sT_{\mathcal{KID}})$ into a dense partition of
$(X,\sT_{\mathcal{KID}})$ of cardinality $\tau$.
It is not necessary to assume that $\sA=\{D^n:n<\omega\}$.

(b) The argument of Theorem~\ref{notkappa'resolv}(b) closely parallels
our proof in \cite{comfhu07}(4.2) that an $\omega$-resolvable, dense
subset $X$ of a space of the form $(D(\kappa))^I$ is necessarily
$\kappa$-resolvable (i.e., is maximally resolvable in case
$\Delta(X)=\kappa$). That theorem,
surprising to the authors, helps to explain the difficulty encountered
over the years by many workers attempting to answer the question of Ceder and
Pearson~\cite{cederpear}: Is every $\omega$-resolvable space maximally
resolvable?

(c) It should be noted that a dense subspace of a space of the form
$(D(\kappa))^I$ need not be $\omega$-resolvable. Indeed in
\cite{comfhu07}(2.3) we show that for every $\kappa\geq\omega$ there is
a dense set $X\subseteq(D(\kappa))^{2^\kappa}$ such that
$|X|=\Delta(X)=\kappa$, no subset of $X$ is resolvable, and every dense
subset of $X$ is open in $X$. See also \cite{asttw}(2.3),
\cite{comfhu03}(5.4) and \cite{juhss}(4.1) for parallel results in the
space $\{0,1\}^{2^\kappa}$.

(d) A propos of (b) above, we note that other criteria sufficient to
ensure maximal resolvability have been established by other authors.
For example, years ago Pytke$'$ev~\cite{pyt83} showed that every
$k$-space, also every space $X$ for which the tightness $t(X)$ satisfies
$t(X)<\Delta(X)$, is maximally resolvable. More recently,
denoting by $\ps(X)$ the smallest successor cardinal such that every
discrete set $S\subseteq X$ satisfies $|S|<\ps(X)$,
Pavlov~\cite{pavlov02} showed that every $T_1$-space such that
$\Delta(X)>\ps(X)$ is maximally resolvable. That theorem was
strengthened in two ways in
\cite{juhss2}: No separation hypothesis on $X$ is required,
and maximal resolvability of $X$ is established assuming only
$\Delta(X)\geq\ps(X)$.
}
\end{remarks}

Our proof of Theorem~\ref{notkappa'resolv} rests on the conventions of
 Section~2, and uses crucially the (strong) hypothesis that $(X,\sT)$ is
maximally resolvable.
That hypothesis can be weakened to the assumption that $(X,\sT)$ is
$\kappa'$-resolvable, with $\kappa'$ regular and
$S(X,\sT)\leq\kappa'\leq|X|=\Delta(X,\sT)=\kappa$, provided that the
equality $2^{\kappa'}=2^{|X|}$ is assumed. Indeed the argument given in
the proof of Theorem~\ref{notkappa'resolv} shows that
$\sU:=\sT_{\mathcal{KID}}$ has properties (a), (b) and (c), with
$\sD=\{D^n_\eta:n<\omega,\eta<\kappa'\}$ a dense partition of $(X,\sT)$,
with $\sI=\{\sI_t:t\in X\times 2^{\kappa'}\}$ a strong small-set-separating,
$\kappa'$-independent family of partitions of $\kappa'$, and with
$\sK=\widetilde{\sM}$ as in Definition~\ref{defnew} with $Z=X$.
We do not know in ZFC whether the hypothesis of
Theorem~\ref{notkappa'resolv} can be weakened.
Specifically we ask:

\begin{question}\label{kappa'resolv?}
{\rm
Let $X=(X, \sT)$ be a crowded Tychonoff space
and let $\kappa'$ be a regular cardinal such that
$S(X,\sT)\leq\kappa'<|X|= \Delta(X,\sT)=\kappa$ and $(X,\sT)$ is
$\tau$-resolvable for each $\tau<\kappa'$. Must there then exist, in ZFC, a
Tychonoff refinement $\sU$ of $\sT$ such that

{\rm (a)} $S(X,\sU)\leq\kappa'$
(perhaps even: $S(X,\sU)=S(X,\sT)$) and $\Delta(X,\sU)=\Delta(X,\sT)=\kappa)$;

{\rm(b)}  $(X, \sU)$ is 
$\tau$-resolvable for each $\tau<\kappa'$; and

{\rm (c)} $(X,\sU)$ is not $\kappa'$-resolvable?
}
\end{question}

Of course, Question~\ref{kappa'resolv?} is of interest only if
$(X,\sT)$ is itself $\kappa'$-resolvable, since otherwise $\sU:=\sT$ would
be as required.

Next we prove item ($i=3$) of the
Abstract for the case $|X|=\Delta(X)$.

\begin{theorem}\label{i=3}
Let $X=(X, \sT)$ be a crowded,
maximally resolvable Tychonoff space 
with $S(X,\sT)\leq|X|= \Delta(X,\sT)=\kappa$.
Then there is a
Tychonoff refinement $\sU$ of $\sT$ such that

{\rm (a)} $S(X,\sU)=S(X,\sT)$ and $\Delta(X,\sU)=\Delta(X,\sT)$;

{\rm(b)}  $(X, \sU)$ is 
maximally resolvable; and

{\rm (c)} $(X,\sU)$ is not extraresolvable.
\end{theorem}

Proof. We invoke the conventions of \ref{defKID} and \ref{shan1}, now
taking $\tau=\kappa$.

Let $\sD=\{D^\gamma_\eta:\eta<\kappa,\gamma<\kappa\}$
be a faithfully indexed dense
partition of $(X,\sT)$, and set
$D^\gamma:=\bigcup_{\eta<\kappa}\,D^\gamma_\eta$
for $\gamma<\kappa$. Let $\sI=\{\sI_t:t\in X\times 2^\kappa\}$ be a
$\kappa$-independent family of partitions
$\sI_t$ of $X$ with the strong small-set-separating property;
for simplicity we take $\kappa_t=2=\{0,1\}$ for each $t\in
X\times2^\kappa$.

Let $\sM=\{M_\xi:\xi<2^\kappa\}=\sP(X)$, and define
$\sK:=\widetilde{\sM}$ as in Definition~\ref{defnew} (taking $Z=X$).
We will show that
$\sU:=\sT_{\mathcal{KID}}$ is as required.

(a) The equality $\Delta(X,\sT_{\mathcal{KID}})=\Delta(X,\sT)$ is given by
Corollary~\ref{delta=delta},
while $S(X,$ $\sT_{\mathcal{KID}})$ $=S(X,\sT)$ is immediate from
Lemma~\ref{S(KID)} (using the regularity of $S(X,\sT)$ and the fact that
$\kappa_t<\omega<\omega^+\leq S(X,\sT)$ for each $t\in Z\times2^\kappa$).

(b) According to Corollary~\ref{delta=delta}(b), the disjoint sets $D^\gamma$
($\gamma<\kappa$) are dense in $(X,\sT_{\mathcal{KID}})$.

(c) Suppose there is a family $\sE$ of dense subsets of
$(X,\sT_{\mathcal{KID}})$, with $|\sE|=\kappa^+$, such that every two
elements of $\sE$ have intersection which is nowhere dense in
$(X,\sT_{\mathcal{KID}})$. We claim that, much as in the proof of
Theorem~\ref{maintheorem}(c), there is $E\in\sE$ such that 

$\int_{(D^\gamma,\sT_{\mathcal{KID}})}(D^\gamma\cap E)=\emptyset$ for each
$\gamma<\kappa$.\hfill(*)

\noindent For if (*) fails then some (fixed)
$\gamma<\kappa$ satisfies

$\int_{(D^\gamma,\sT_{\mathcal{KID}})}(D^\gamma\cap E)\neq\emptyset$ for
$\kappa^+$-many $E\in\sE$,\\
\noindent and then since
$S(D^\gamma,\sT_{\mathcal{KID}})=S(X,\sT_{\mathcal
KID})=S(X,\sT)\leq\kappa$
there are distinct $E,E'\in\sE$ such that

$\emptyset\neq[\int_{(D^\gamma,\sT_{\mathcal{KID}})}(D^\gamma\cap E)]\cap
[\int_{(D^\gamma,\sT_{\mathcal{KID}})}(D^\gamma\cap E')]=
\int_{(D^\gamma,\sT_{\mathcal{KID}})}(D^\gamma\cap E\cap E')$.\\
\noindent Then with $\sT_{\mathcal{KID}}$-open $U\subseteq X$ chosen so
that
$D^\gamma\cap U=\int_{(D^\gamma,\sT_{\mathcal{KID}})}(D^\gamma\cap E\cap
E')$\\
\noindent we have

$\emptyset\neq U\subseteq\cl_{(X,\sT_{\mathcal{KID}})}U
=\cl_{(X,\sT_{\mathcal{KID}})}(D^\gamma\cap U)
=\cl_{(X,\sT_{\mathcal{KID}})}\int_{(X,\sT_{\mathcal{KID}})}(D^\gamma\cap
E\cap E')
\subseteq\cl_{(X,\sT_{\mathcal{KID}})}(E\cap E')$,\\
\noindent contrary to the fact that $E\cap E'$ is nowhere dense in
$(X,\sT_{\mathcal{KID}})$. Thus (*) is established.

Then, choosing $E\in\sE$ as in (*), we have from
Theorem~\ref{heredirr}(b) (applied to the set $M_{\overline{\xi}}=E$)
that $E\in\sK=\widetilde{\sM}$, so by Lemma~\ref{KIDprops}((a) and (b))
the set $E$ is closed and discrete in the crowded space $(X,\sT_{\mathcal
KID})$. This contradicts the density of $E$ in
$(X,\sT_{\mathcal{KID}})$.~$\square$

We turn next to establishing items  ($i=4$) and ($i=5$) of the Abstract
for the case $|X|=\Delta(X)$.
As expected, refinements of
the form $\sU=\sT_{\mathcal{KID}}$ play a central role; it
is necessary only to tailor in each case the specifics of the
families $\sK$, $\sI$, and $\sD$ to the task at hand.
But in Theorem~\ref{i=5} the process is iterated:
a first expansion
$\sT'\supseteq\sT$ satisfies $\nwd(X,\sT')=\kappa$, a second expansion
$\sT''\supseteq\sT'$ is maximally resolvable but not extraresolvable,
and a final expansion (of the form $\sT''^{\mathcal{ID}}$, not
$\sT''_{\mathcal{KID}}$) has all required properties.

For the proofs of (the case $|X|=\Delta(X)$ of) items ($i=4$)
and ($i=5$) of the Abstract, we need two preliminary lemmas.
A weak version of Lemma~\ref{S_er} is proved in
our work \cite{comfhu04}(3.9). A strictly combinatorial
proof exists, but it is lengthy; we give instead an argument which uses
the topological constructions already at our disposal.

\begin{lemma}\label{S_er}
Let $\tau\geq\omega$. There exist families
$\sA=\{\sA_\xi:\xi<2^\tau\}$ and $\sS_{er}\subseteq\sP(\tau)$ such
that

{\rm (i)} $\sA$ is a $\tau$-independent family of
partitions of $\tau$ with the strong small-set-separating property,
with each $\sA_\xi\in\sA$ of the form $\sA_\xi=\{A^0_\xi,A^1_\xi\}$;

{\rm (ii)} $|\sS_{er}|=2^\tau$;

{\rm (iii)}
if $n<\omega$ and $S, S_1, S_2,\ldots S_n$ are distinct elements of
$\sS_{er}$ and
$A=\bigcap_{\xi\in F}\,A^{f(\xi)}_\xi$
with $F\in[2^\tau]^{<\omega}$ and $f\in\{0,1\}^F$,
then $|A\cap(S\backslash(S_1\cup S_2\cup\ldots\cup S_n))|=\tau$; and

{\rm (iv)} 
if $S,S'\in\sS_{er}$ with $S\neq S'$ then
{\rm (a)}~for each $x\in\tau\backslash(S\cap S')$ there are infinitely
many $\xi<2^\tau$ such that
$x\in A^1_\xi$ and $S\cap S'\subseteq A^0_\xi$; and {\rm (b)}~for each
$x\in S\cap S'$ there are infinitely many
$\xi<2^\tau$ such that $(S\cap S')\cap A^1_\xi=\{x\}$.

\end{lemma}

Proof. Let $\sJ\cup\sL\cup\{\sD\}$ be a $\tau$-independent family of
partitions of $\tau$, 
where $\sJ=\{\sJ_\xi:\xi<2^\tau\}$
is chosen (as in Theorem~\ref{nwdX=kappa}) so that the
space

$Y=(Y,\sT):=e_\sJ[\tau]\subseteq K:=\{0,1\}^\sJ=\{0,1\}^{2^\tau}$\\
\noindent has properties (a), (b), (c) and (d) of Theorem~\ref{nwdX=kappa}.
We take $|\sJ|=|\sL|=2^\tau$, say
$\sJ=\{\sJ_\xi:\xi<2^\tau\}$ and
$\sL=\{\sL_\zeta:\zeta<2^\tau\}$, and we take each
$\sJ_\xi\in\sJ$ of the form $\sJ_\xi=\{J^0_\xi,J^1_\xi\}$ and each
$\sL_\zeta\in\sL$ of the form
$\sL_\zeta=\{L^0_\zeta,L^1_\zeta\}$.

We write $\sD=\{D^\gamma_\eta:\gamma<\tau,\eta<\tau\}$.

The families $\sA$ and $\sS_{er}$ will be defined with the
help of a suitable expansion $\sT_{\mathcal KID}$ of $\sT$.

The family
$\sD$ has already been defined, and for $\sI$ we choose an arbitrary
$\tau$-independent
family $\sI=\{\sI_t:t\in Y\times2^\tau\}$
of partitions of $\tau$ with the strong small-set-separating property,
say with each $\sI_t$ of the form $\sI_t=\{I^0_t,I^1_t\}$.
For $\sK$, first let

$\sK':=\{\bigcap_{\zeta\in F}\,L^0_\zeta:|F|>1,
F\in[2^\tau]^{<\omega}\}$\\
\noindent and let $\sK$ be the set of sets of the form
$\bigcup_{i<n}\,K_i'$ with $n<\omega$, $K_i'\in\sK'$.
We write $\sK=\{K_\xi:\xi<2^\tau\}$, the indexing chosen so that
each $K\in\sK$ is listed infinitely often.

With these definitions, writing as usual
$D^\gamma:=\bigcup_{\eta<\tau}\,D^\gamma_\eta$ for $\gamma<\tau$,
conditions
(1), (2), (3) and (4) of \ref{shan1} are clearly satisfied.
To verify (5), fix $K\in\sK$ and $\gamma<\tau$; we show that
$\int_{(D^\gamma,\sT^{\mathcal{ID}})}\,(K\cap D^\gamma)=\emptyset$.

There are $n<\omega$ and $F_i\in[2^\tau]^{<\omega}$ with $|F_i|>1$
such that $K=\bigcup_{i<n}(\bigcap_{\zeta\in F_i}\,L^0_\zeta)$. With
$F:=\bigcup_{i<n}\,F_i$ and $B:=\bigcap_{\zeta\in F}\,L^1_\zeta$ we have
(since $\sL\cup\{\sD\}$ is $\tau$-independent)
for each $\eta<\tau$ that $|B\cap D^\gamma_\eta|=\tau$, so
$B\cap D^\gamma$ meets each set of the form
$H^\alpha_t=X(I^\alpha_t)$ with $I^\alpha_t\in\sI_t\in\sI$.
Thus $B\cap D^\gamma$
is dense in $(D^\gamma,\sT^{\mathcal{ID}})$, and
from $B\cap K=\emptyset$ it then follows that 
$\int_{(D^\gamma,\sT^{\mathcal{ID}})}\,(K\cap D^\gamma)=\emptyset$. Thus
(5) is proved.

Now with $\sW=\{\sW_t:t\in Y\times2^\tau\}$
defined (using $\sI$ and $\sK$)
as in Definition~\ref{defKID} we set
$\sA:=\sJ\cup\sW$,
and $\sS_{er}:=\{L^0_\zeta:\zeta<2^\tau\}$. It is clear for distinct
$S,S'\in\sS_{er}$ that $S\cap S'\in\sK'\subseteq\sK$.

Each $\sJ_\xi\in\sJ$, and
each $\sW_t\in\sW$, is a partition of $\tau$. Since $\sJ$
has the strong small-set-separating property, and $\sJ\subseteq\sA$,
also $\sA$ has the strong small-set-separating property. Thus
to prove (iii) and to complete the proof of (i) it suffices to show:
For each triple $(J,W,L)$, with

$J=\bigcap_{\xi\in F_0}\,J_\xi^{f_0(\xi)}$ with
$F_0\in[2^\tau]^{<\omega}$, $f_0\in\{0,1\}^{F_0}$,

$W=\bigcap_{t\in F_1}\,W_t^{f_1(t)}$ with
$F_1\in[Y\times2^\tau]^{<\omega}$, $f_1\in\{0,1\}^{F_1}$, and

$L=L_{\overline{\zeta}}^0\backslash\bigcup_{i<n}\,L_{\zeta_i}^0
=L^0_{\overline{\zeta}}\cap\bigcap_{i<n}\, L^1_{\zeta_i}$
with distinct $\overline{\zeta},\zeta_i<2^\tau$,\\
\noindent that $|J\cap W\cap L|=\tau$.

To do that, take $|F_1|=m$, say $F_1=\{t_j=(x_j,\xi_j):j<m\}$,
and note
with $K_{\xi_j}=\bigcup_{i<n_j}(\bigcap_{\zeta\in F_{i,j}}\,L^0_\zeta)$
that $L\backslash K_{\xi_j}$ contains the set
$C:=L\cap\bigcap\{L^1_\zeta:\zeta\in\bigcup_{i<n_j} F_{i,j}\}$. Thus
$L\cap W\supseteq C\cap W=C\cap H$, where
$H=\bigcap_{t\in F_1}\,H_t^{f_1(t)}$. Since
$\sJ\cup\sL\cup\{\sD\}$ is $\tau$-independent,
and each $H\in\sH_t\subseteq\sT^{\mathcal{ID}}$ is the union of sets
in $\sD$, the family $\sJ\cup\sH\cup\sL$ is also $\tau$-independent.
Now $J$ is a Boolean combination of sets from $\sJ$,
$H$ is a Boolean combination of sets from $\sH$, and
$C$ is a Boolean combination of sets from $\sL$,
so from $J\cap W\cap L\supseteq J\cap W\cap C=J\cap H\cap C$ then follows
$|J\cap W\cap L|=\tau$, as required.

For (iv), let $S,S'\in\sS_{er}$ with $S\neq S'$ and
fix $x\in Y$.
Then for each of the (infintely many) $\xi<2^\tau$ such that $S\cap
S'=K_\xi$ we have, taking $t=(x,\xi)$: If $x\notin S\cap S'$ then
$x\in W^1_t$ and $S\cap S'\subseteq W^0_t$
with $W^0_t,W^1_t\in\sW_t\subseteq\sA$, while if $x\in S\cap S'$ then
$(S\cap S')\cap W^1_t=\{x\}$ with $W^1_t\in\sW_t\subseteq\sA$.
Then to achieve (iv) in the form stated, it is enough to re-index $\sA$
in the form $\sA=\{\sA_\xi:\xi<2^\tau\}$~$\square$

\begin{theorem}\label{usingS_er}
Let $\tau\leq\kappa$ and let $(X,\sT)$ be a crowded, $\tau$-resolvable
Tychonoff space such that $S(X,\sT)\leq|X|=\Delta(X,\sT)=\kappa$. Then
there is a Tychonoff expansion $\sU$ of $\sT$ such that

{\rm (a)} $S(X,\sU)=S(X,\sT)$ and $\Delta(X,\sU)=\Delta(X,\sT)$;

{\rm (b)} $(X,\sU)$ is $\tau$-resolvable; and

{\rm (c)} $(X,\sU)$ is $2^\tau$-extraresolvable.
\end{theorem}

Proof. 
Let $\sA$ and $\sS_{er}$
be families 
given by Lemma~\ref{S_er}.
Ignoring the indexing there,
we choose $Z\subseteq X$ with $|Z|=1$
and we write
$\sA=\{\sA_t:t\in Z\times2^\tau\}$, with each $\sA_t=\{A^0_t,A^1_t\}$.
Let $\sK=\{K_\xi:\xi<2^\tau\}$ with each
$K_\xi=\emptyset$, and
let $\sD:=\{D_\eta^\gamma: \eta<\tau,\gamma<\tau\}$ be a partition
of $X$
witnessing the $\tau$-resolvability of 
$(X, \sT)$. We show that the expansion
$\sU:=\sT^{\mathcal{AD}}$ satisfies
conditions (a), (b), and (c).

(a) That $\Delta(X,\sU)=\Delta(X,\sT)$ is given by
Corollary~\ref{delta=delta}(c), while
$S(X,\sU)=S(X,\sT)$ is immediate from
Lemma~\ref{S(KID)} (using the regularity of $S(X,\sT)$ and the fact that
$\kappa_t=2<\omega<\omega^+\leq S(X,\sT)$ for each $t\in Z\times2^\tau$).
Thus (a) holds.

(b) As usual, Lemma~\ref{KIDprops}(c)
shows that $\{D^\gamma:\gamma<\tau\}$ is a dense partition of $(X,\sU)$.

(c) It suffices to show that

~~~(i) if $S\in\sS_{er}$ then $X(S)$ is dense in $(X,\sU)$; and

~~~(ii) if $S,S'$ are distinct elements of $\sS_{er}$ then $X(S\cap S')$
is closed and nowhere dense in $(X,\sU)$.

For (i), given $\emptyset\neq U\in\sT$ and
$H=\bigcap_{t\in F}\,H^{f(t)}_t$ with $F\in[Z\times2^\tau]^{<\omega}$
and $f\in\{0,1\}^F$, we must show $X(S)\cap(U\cap H)\neq\emptyset$.
Set $A=\bigcap_{t\in F}\,A^{f(t)}_t$, so that $H=X(A)$. Then $A\cap
S\neq\emptyset$
(indeed $|A\cap S|=\tau$ by Lemma~\ref{S_er}(iii)) so there are
$\tau$-many pairs $(\gamma,\eta)$ such that $D^\gamma_\eta\subseteq
X(A)\cap X(S)$. Each such $D^\gamma_\eta$ meets $U$, so

$|X(S)\cap(H\cap U)|=|X(S)\cap X(A)\cap U|=\tau$.

For (ii), let $p\in X\backslash X(S\cap S')$, say $p\in D^\gamma_\eta$,
and using
Lemma~\ref{S_er}(iv) choose
$t=(x,\xi)\in Z\times2^\tau$ such that $\eta\in A^1_t$ and $S\cap S'\subseteq
A^0_t$. Then $p\in X(A^1_t)=H^1_t$ and $X(S\cap S')\subseteq
X(A^0_t)=H^0_t$, so $H^1_t$ is a $\sU$-open neighborhood of
$p$ disjoint from $X(S\cap  S')$. Thus $X(S\cap S')$ is closed in
$(X,\sU)$.

Given $\emptyset\neq U\in\sT$ and $H=\bigcap_{t\in F}\,H^{f(t)}_t$ as
in (a) set $A:=\bigcap_{t\in F}\,A^{f(t)}_t$, so that $H=X(A)$, and note
from Lemma~\ref{S_er}(iii) that $|A\cap(S'\backslash S)|=\tau$. Then

$U\cap X(A)\cap X(S'\backslash S)=U\cap H\cap X(S'\backslash
S)\neq\emptyset$,\\
\noindent so $X(S'\backslash S)=X(S')\backslash X(S)$ is dense in
$(X,\sU)$. {\it A fortiori} $X\backslash X(S)$ is dense in $(X,\sU)$,
so the closed set $X(S)$ is nowhere dense in $(X,\sU)$.~$\square$

Now we are ready to prove the case $|X|=\Delta(X)$ of
items ($i=4$) and ($i=5$) of the Abstract.

\begin{theorem}\label{i=4}
Let $X=(X, \sT)$ be a crowded,
maximally resolvable Tychonoff space 
with $S(X,\sT)\leq|X|= \Delta(X,\sT)=\kappa$.
Then there is a
Tychonoff refinement $\sU$ of $\sT$ such that

{\rm (a)} $S(X,\sU)=S(X,\sT)$ and $\Delta(X,\sU)=\Delta(X,\sT)$;

{\rm (b)} $(X,\sU)$ is extraresolvable; and

{\rm (c)} $(X,\sU)$ is not maximally resolvable.
\end{theorem}

Proof. The topology $\sU$ will be of the form $\sU=\sT_{\mathcal{KAD}}$.
We first define the families $\sK$, $\sA$ and $\sD$.

Let $\sD=\{D_\eta:\eta<\kappa\}$ be a dense partition of $(X,\sT)$ which
witnesses the maximal resolvability of $(X,\sT)$. (Note. To match the
notation used throughout Section~2, more formally we take
$\tau=1=\{0\}$ and $D_\eta=D_\eta^0$ in Notation~\ref{X(S)}; then
$X(S)=\bigcup_{\eta\in S}\,D_\eta$ for $S\subseteq\kappa$.)

Let $\sA=\{\sA_\xi:\xi<2^\kappa\}$ with $\sA_\xi=\{A_\xi^0, A_\xi^1\}$
and
$\sS_{er}\subseteq\sP(\kappa)$ as given in Lemma~\ref{S_er}, and re-index
$\sA$ in the form 
$\sA = \{\sA_t: t\in X\times2^\kappa\}$. We partition the set $2^\kappa$
in the form
$2^\kappa=T_0\cup T_1$ with $|T_0| = |T_1| = 2^\kappa$. We assume
without  loss 
of generality that the families $\{\sA_t: t\in X\times T_1\}$
and $\sS_{er}$ satisfy
conditions (i) through (iv) of Lemma~\ref{S_er}.

The definition of the family $\sK$ parallels the construction in
Definition~\ref{defnew}, but with modifications. Specifically:

Let $\sM=\{M_\xi:\xi<2^\kappa\}=\sP(X)$ with $M_0=\emptyset$ and define
$\widetilde{M}=\{\widetilde{M_\xi}:\xi<2^\kappa\}$ as follows.

$\widetilde{M_0}=\emptyset$; and\\
\noindent if $0<\xi<2^\kappa$ and $\widetilde{M_\eta}$ has been defined
for all $\eta<\xi$ then

$\widetilde{M_\xi}=M_\xi$ if each set of the form

$(M_\xi\cup\widetilde{M_{\eta_0}}\cup
\widetilde{M_{\eta_1}}\cup\cdots\cup
\widetilde{M_{\eta_n}})\cap X(S)$
($n<\omega$, $\eta_i<\xi$, $S\in\sS_{er}$)\\
\noindent has nonempty interior in the space $(X(S),\sT^{\mathcal{AD}})$

$=\emptyset$ otherwise.

Then with $T_0,T_1\subseteq2^\kappa$ as above, we write
$\sK=\sK_0\cup\sK_1$ with $\sK_i=\{K_\xi:\xi\in T_i\}$; we arrange that
$\{K_\xi:\xi\in T_0\}$ is a faithful indexing of 
$\widetilde{\sM}$, and $K_\xi=\emptyset$ for each $\xi\in T_1$.

We claim that the topology $\sU = \sT_{\sK\sA\sD}$ is as required.

We verify conditions (a), (b) and (c). Indeed as to (c) we will show
that $(X,\sU)$ is not even $S(X,\sT)$-resolvable.

(a) From Corollary~\ref{delta=delta} and Lemma~\ref{S(KID)} we have

$$\kappa=\Delta(X,\sT)=\Delta(X,\sT_{\mathcal{KAD}})=
\Delta(X,\sU)~{\mathrm{and}}~\kappa=S(X,\sT)
=S(X,\sT_{\mathcal{KAD}})=S(X,\sU).$$

(b) It is enough to show that

(i) if $S\in\sS_{er}$ then $X(S)$ is dense in $(X,\sU)$; and

(ii) if $S$ and $S'$ are distinct elements of $\sS_{er}$, then $X(S)\cap
X(S')$ is (closed and) nowhere dense in $(X,\sU)$.

For (i), we must show that if $\emptyset\neq U\in\sT$ and
$W=\bigcap_{t\in F}\,H_t^{f(t)}\backslash K_t\in\sT_{\mathcal{KAD}}$
with $F\in[X\times2^\kappa]^{<\omega}$ and $f\in\{0,1\}^F$,
 then $X(S)\cap(U\cap
W)\neq\emptyset$. For that, set $A:=\bigcap_{t\in F}\,A_t^{f(t)}$,
so that $H:=\bigcap_{t\in F}\,H_t^{f(t)}=X(A)$ and
$W=X(A)\backslash K$ with $K=\bigcup\{K_t: t\in F\}$. Since
$\int_{(X(S),\sT^{\mathcal{AD}})}(K\cap X(S))=\emptyset$, the set
$X(S)\backslash K$ is dense in $(X(S),\sT^{\mathcal{AD}})$. 
From Lemma~\ref{S_er}(iii) we have $A\cap S\neq\emptyset$, so
$X(A)\cap U\cap X(S) = X(A\cap S)\cap U\neq\emptyset$. Since  
$\emptyset\neq X(A)\cap U\in\sT^{\sA\sD}$ 
and $K$ is closed and discrete in $\sU$,
we have 
$(X(S)\setminus K)\cap X(A)\cap U\neq\emptyset$ and therefore
$X(S)\cap W\cap U\neq\emptyset$, as required.

Before verifying (ii), we show this
for later use.

each set $\bigcup_{\eta\in G}\,D_\eta$ with
$G\in[\kappa]^{<\omega}$ is in $\sK$.\hfill (1)\\
\noindent If that fails, there are $U\in\sT$, $K\in\sK$, $S\in\sS_{er}$, and
$H\in\sT^{\mathcal{AD}}$ such that

$U\cap H\cap X(S)\subseteq[(\bigcup_{\eta\in G}\,D_\eta)\cup K]\cap
X(S)$.\\ 
\noindent Here $H=X(A)$ with $A=\bigcup_{t\in F}\,A_t^{f(t)}$.
Since $\{\sA_t: t\in X\times T_1\}$ has the strong
small-set-separating property,
for each $\eta\in G$ there are infinitely many
indices $v_\eta$ such that ($\sA_{v_\eta}\in\sA$ and) $\sA_{v_\eta}$
separates $\{\eta\}$ and $\emptyset$. For each $\eta\in G$ we choose
such $v_\eta$ such that $v_\eta\notin F$, so that
$A\cap \bigcap_{\eta\in G}\,A^1_{v_\eta}\neq\emptyset$. We set 
$H':= \bigcap_{\eta\in G}\,X(A^1_{v_\eta})$, so that
$\emptyset\neq U\cap H\cap H'\in\sT^{\sA\sD}$. Then

$[(\bigcup_{\eta\in G}\,D_\eta)\cup K]\cap X(S)\supseteq 
U\cap H\cap X(S)\supseteq
U\cap H\cap H'\cap X(S)\neq\emptyset$.\\
\noindent Here 
$U\cap H\cap H'\cap X(S)\neq\emptyset$ because 
$X(S)$ is dense in $\sU$ (see (b)(i)) and 
$H'$ differs from certain $W'\in\sU$ by a 
$\sU$-closed, $\sU$-discrete set $K\in\sK$.

Since
$(\bigcap_{\eta\in G}\,A^1_{v_\eta})\cap(\bigcap_{\eta\in F}\,D_\eta)=
\emptyset$, we then have

$K\cap X(S)\supseteq U\cap H\cap H'\cap X(S)\neq\emptyset$,

\noindent contradicting the condition $K\in\sK$. Thus (1) is shown.

Now for (ii), let $x\in X\backslash X(S\cap S')$, say with $x\in D_\eta$,
and using Lemma~\ref{S_er}(iv)(a) choose $u\in X\times T_1$ such that
$\eta\in A^1_u$ and $S\cap S'\subseteq A^0_u$. Then $x\in
X(A^1_u)=H^1_u = W^1_u$ and $X(S)\cap X(S')\subseteq X(A^0_u)=H^0_u=W^0_u$ 
(since $K_\xi=\emptyset$ for $\xi\in T_1$ in Definition~\ref{defKID}),
so $W^0_u$ is
a neighborhood in $(X,\sU)$ of $x$ which is disjoint from $X(S\cap S')$.
Thus $X(S\cap S')$ is closed in $(X,\sU)$.

To see that the closed set $X(S\cap S')$ is nowhere dense in $(X,\sU)$,
suppose (taking notation as above) that there are nonempty $U\in\sT$ and
$W=X(A)\backslash K\in\sU$ with 
$A=\bigcap_{t\in F}\,A_t^{f(t)}$ such that
$U\cap W\subseteq X(S\cap S')$. Fix $\eta\in S\cap S'$
and use Lemma~\ref{S_er}(iv)(b) to find
$u\in X\times T_1$ such that $u\notin F$ 
and $(S\cap S')\cap A^1_u=\{\eta\}$. Then
$X(A^1_u)\cap X(S\cap S')=D_\eta$, and the condition $u\notin F$ implies 
$\emptyset\neq X(A)\cap X(A^1_u)\in\sT^{\mathcal{AD}}$, which further implies 
$\emptyset\neq U\cap W\cap X(A^1_u)\in\sU$.  Hence 
$U\cap W\cap X(A^1_u)\subseteq 
U\cap W\subseteq X(S\cap S')$ and from (1) we have

$\emptyset\neq U\cap W\cap X(A^1_u)\subseteq X(A^1_u)\cap X(S\cap S')
\subseteq D_\eta\cup K\in\sK$.\\ 
\noindent But from Lemma~\ref{KIDprops} the space
$D_\eta\cup K$
is closed and discrete in $(X,\sU)$, a contradiction. The proof of (b)
is complete.
 
(c) Here we show more, namely that $(X,\sU)$ is not even
$S(X,\sT)$-resolvable. Arguing much as in Theorem~\ref{heredirr}(b), we
first show this:

if $\xi<2^\kappa$ and $\int_{(X(S),\sU)}(M_\xi\cap X(S))=\emptyset$
for all $S\in\sS_{er}$, then $\widetilde{M_\xi}=M_\xi\in\sK$.\hfill (2)

\noindent For that, we must show for fixed $S\in\sS_{er}$ and fixed $K\in\sK$
that

$\int_{(X(S),\sT^{\mathcal{AD}})}[(M_\xi\cup K)\cap X(S)]=\emptyset$.

To see that, let $\emptyset\neq U\in\sT$ and $W=\bigcap_{t\in
F}\,W_t^{f(t)}\in\sU$.
Since $X(S)\setminus M_\xi$ is dense in $(X(S),\sU)$ and $\emptyset\neq U\cap
W\in\sU$, we have that
$Y:=(X(S)\backslash M_\xi)\cap(U\cap W)$ is dense in
$(X(S)\cap(U\cap W),\sU)$. Thus $Y$ is crowded, so since $W$
differs from $H:=\bigcap_{t\in F}\,H_t^{f(t)}\in\sT^{\mathcal{AD}}$
by a set $K'\in\sK$
we have from Lemma~\ref{discrete} that $Y\backslash(K'\cup
K)$ remains dense in $(X(S)\cap(U\cap W),\sU)$, hence dense in
$(X(S)\cap(U\cap H),\sT^{\mathcal{AD}})$. Thus

$\int_{(X(S),\sT^{\mathcal{AD}})}[(M_\xi\cup K)\cap X(S)]=\emptyset$,\\
\noindent as required, and (2) is proved.

To complete the proof of (c) we argue by contradiction, supposing that
$\{E_\eta:\eta<S(X,\sT)\}$ is a pairwise disjoint family of dense
subsets of $(X,\sU)$. For each $\eta<S(X,\sT)$ there is
$S_\eta\in\sS_{er}$ such that $\int_{(X(S_\eta),\sU)}(E_\eta\cap
X(S_\eta))\neq\emptyset$, so there are nonempty $U_\eta\in\sT$ and
$W_\eta=X(A_\eta)\backslash K_\eta\in\sU$ with 
$A_\eta=\bigcap_{t\in F_\eta}\,A_t^{f_\eta(t)}$
such that

$\emptyset\neq U_\eta\cap W_\eta\cap X(S_\eta)\subseteq
E_\eta\cap X(S_\eta)$.

\noindent For notational simplicity set
$V_\eta:=U_\eta\cap W_\eta\cap X(S_\eta)$
for $\eta<S(X,\sT)$. Then

$\{V_\eta:\eta<S(X,\sT\}$ is a pairwise
disjoint family,\hfill (3)\\
\noindent since $V_\eta\subseteq E_\eta$ and
$\{E_\eta:\eta<S(X,\sT\}$ is pairwise disjoint.

Now recall, using the notation
$\sJ$, $\sW$, $\sH_t$, $\sL$, $\sD$, $J$, $W$, $H$, $K$, $F$
and $L$ as in (the proof of)
Lemma~\ref{S_er}, that
each of the present sets
$S_\eta$ is of the form $L_\zeta^0\in\sL_\zeta\in\sL$ for
some $\zeta<2^\kappa$, and
$A_\eta$
is a Boolean combination of sets from $\sJ\cup\sW$, say
$A_\eta=J\cap W$ where $J, W$ are as in the proof of Lemma~\ref{S_er}.
Write $W = H\setminus K$ with $H$ as in Lemma~\ref{S_er} and 
with $K$  of the form
$K=\bigcup_{i<n}(\bigcap_{\zeta\in F_i}\,L^0_\zeta)$ with
$1<|F_i|<\omega$. For each $i<n$ choose
$\zeta_i\in F_i$ such that $L_{\zeta_i}^0\neq S_\eta$. Then 
$L:=\bigcap_{i<n} L_{\zeta_i}^1$
satisfies $L\cap K=\emptyset$ and $L\cap S_\eta\neq\emptyset$. 
  
Since $\sJ\cup\sL\cup\{\sD\}$ is an independent family and
the elements of each partition in $\sH$ are unions of 
some dense sets in $\sD$, the family
$\sJ\cup\sH\cup\sL$ is also an independent family. Since
$J$ is a Boolean combination of sets from $\sJ$, 
$H$ is a Boolean combination of sets from $\sH$,
and $L\cap S_\eta$ is a Boolean combination of sets from $\sL$,
we have $J\cap H\cap L\cap S_\eta\neq \emptyset$.
Since $W = H\setminus K$ and $L\cap K=\emptyset$, we
have $H\cap L\subseteq W$ and 

$$\emptyset\neq J\cap H\cap L\cap S_\eta = J\cap (H\cap L)\cap S_\eta
\subseteq J \cap  W\cap S_\eta = A_\eta\cap S_\eta. 
$$

This argument shows that for each $\eta<S(X,\sT)$ there is a Boolean
combination $N_\eta$ 
of sets from the independent family
$\sJ\cup\sH\cup\sL$, of the form
$N_\eta=P\cap H\cap L\cap S_\eta$, 
such that

$\emptyset\neq N_\eta\subseteq A_\eta\cap S_\eta$.\hfill (4)\\
\noindent For simplicity write $\sB:=\sJ\cup\sH\cup\sL=\{\sB_t:t\in T\}$
with $|T|=2^\kappa$ and write each $N_\eta$ in the form
$N_\eta=\bigcap_{t\in F_\eta}\,B_t^{i_\eta(t)}$ with
$F_\eta\in[T]^{<\omega}$, $i_\eta\in\{0,1\}^{F_\eta}$. Since $S(X,\sT)$
is a regular cardinal there are, by the Erd{\H{o}}s-Rado theorem on
quasi-disjoint sets \cite{comfneg74}, \cite{comfneg82} (the
``$\Delta$-system Lemma"~\cite{juhasz80}) a (finite) set $F$ and
$Q\subseteq S(X,\sT)$ with $|Q|=S(X,\sT)$ such that $F_\eta\cap
F_{\eta'}=F$ whenever $\eta,\eta'\in Q$, $\eta\neq\eta'$. We assume
without loss of generality that $F\neq\emptyset$ and that
$i_\eta(t)=i_{\eta'}(t)\in\{0,1\}$ for all $\eta,\eta'\in Q$,
$t\in F$. Then

$\emptyset\neq N_\eta\cap N_{\eta'}\subseteq(A_\eta\cap
S_\eta)\cap(A_{\eta'}\cap S_{\eta'})$\\
\noindent for distinct $\eta,\eta'\in Q$.

Since $\emptyset\neq U_\eta\in\sT$ for each $\eta\in Q$, there are
distinct $\eta_0,\eta_1\in Q$ (henceforth fixed) such that
$U_{\eta_0}\cap U_{\eta_1}\neq\emptyset$.

Now $\emptyset\neq U_{\eta_k}\cap X(A_{\eta_k})\in\sT^{\mathcal{AD}}$,
and $X(S_{\eta_k})$ is dense in $(X,\sU)$, and $K_{\eta_k}$ is closed and
nowhere dense in $(X,\sU)$, so from $U_{\eta_0}\cap
U_{\eta_1}\neq\emptyset$ follows

$[U_{\eta_0}\cap X(A_{\eta_0})\cap X(S_{\eta_0})\backslash K_{\eta_0}]\cap
[U_{\eta_1}\cap X(A_{\eta_1})\cap X(S_{\eta_1})\backslash K_{\eta_1}]\neq
\emptyset$,\\
\noindent that is:

$V_{\eta_0}\cap V_{\eta_1}=
[U_{\eta_0}\cap W_{\eta_0}\cap X_{\eta_0}]\cap
[U_{\eta_1}\cap W_{\eta_1}\cap X_{\eta_1}]\neq\emptyset$,\\
\noindent which contradicts (3).~$\square$

\begin{theorem}\label{i=5}
Let $X=(X, \sT)$ be a crowded,
maximally resolvable Tychonoff space 
with $S(X,\sT)\leq|X|= \Delta(X,\sT)=\kappa$.
Then there is a Tychonoff refinement $\sU$ of $\sT$ such that

{\rm (a)} $S(X,\sU)=S(X,\sT)$ and $\Delta(X,\sU)=\Delta(X,\sT)$;

{\rm (b)} $(X,\sU)$ is maximally resolvable;

{\rm (c)} $(X,\sU)$ is extraresolvable; and

{\rm (d)} $(X,\sU)$ is not strongly extraresolvable.
\end{theorem}

Proof. We expand in three steps with (modified) $\mathcal{KID}$-like
expansions $\sT\subseteq\sT'\subseteq\sT''\subseteq\sU$. Here are the
details.

Step 1. Let $\sD_1=\{D^\gamma_\eta:\eta<\kappa,\gamma<\kappa\}$ be a
partition of $X$ into $\sT$-dense subsets, let
$\sI_1=\{\sI_t:t\in X\times2^\kappa\}$ be a $\kappa$-independent family of
partitions of $X$ with the strong small-set-separating property with
each $\kappa_t<S(X,\sT)$, and let
$\sK_1:=\{K_\xi:\xi<2^\kappa\}=[X]^{<\kappa}$ (with repetitions
permitted in the indexing of $\sK_1$). Clearly conditions (1), (2), (3)
and (4) of \ref{shan1} are satisfied (with $Z=X$).
To see that (5) also is satisfied,
fix nonempty $U\in\sT$ and $H=X(I)$ with 
$F\in[X\times2^\kappa]^{<\omega}$, $f\in\Pi_{t\in F}\,\kappa_t$,
and $I:=\bigcap_{t\in F}\,I^{f(t)}_t$.  We have $|I|=\kappa$,
so $H\cap D^\gamma\supseteq D^\gamma_\eta$ for
$\kappa$-many $\eta<\kappa$, each dense in $(X,\sT)$, so $|D\cap U\cap
H|=\kappa$. Thus $D^\gamma\cap U\cap H\subseteq K_\xi\cap D^\gamma$ is
impossible, so (5) holds. It follows that
$\sT':=\sT_{\sK_1\sI_1\sD_1}$ has the properties given in
Lemma~\ref{KIDprops}, in particular each $K_\xi\in\sK_1=[X]^{<\kappa}$ is
closed and discrete in $(X,\sT')$, hence nowhere dense, so
$\nwd(X,\sT')=\kappa$.

Step 2. Apply Theorem~\ref{i=3} (to the space $(X,\sT')$) to find an
expansion $\sT''\supseteq\sT'$ (with $\sT''$ of the form
$\sT''=\sT'_{\mathcal{KID}}$) such that $S(X,\sT'')=S(X,\sT')$,
$\Delta(X,\sT'')=\Delta(X,\sT')$, and $(X,\sT'')$ is maximally resolvable
but not extraresolvable.

Step 3. By Theorem~\ref{usingS_er} with $\tau=\kappa$ and with $\sT''$
replacing $\sT$ there, there is an expansion $\sU\supseteq\sT''$ (with
$\sU$ of the form $\sT''^{\mathcal{AD}}$) such that

$S(X,\sU)=S(X,\sT'')=S(X,\sT)$ and
$\Delta(X,\sU)=\Delta(X,\sT'')=\Delta(X,\sT)$\\
\noindent and such that $(X,\sU)$ is
maximally resolvable and $2^\kappa$-extraresolvable. Furthermore
each set $K\in[X]^{<\kappa}$ is closed and discrete
in $(X,\sT')$, hence in
$(X,\sT'')$, so
$\nwd(X,\sT'')=\kappa$. Thus any family $\sE$ (with
$|\sE|>\Delta(X,\sU)=\Delta(X,\sT')$) witnessing the strong
extraresolvability
of $(X,\sU)$ would witness the strong extraresolvability of $(X,\sT'')$,
contrary to the fact that $(X,\sT'')$ is not (even) extraresolvable.~$\square$

\section{The General Case}

The five principal results proved in Section 3 require, in addition to
the essential overarching hypothesis $S(X,\sT)\leq\Delta(X,\sT)$, also the
artificial condition $|X|=\Delta(X,\sT)$. Since for each of those
five results it is essentially
the same argument which allows us to pass from the special
case ($|X|=\Delta(X,\sT)$) to the unrestricted case ($|X|$ is
arbitrary),
we corral all five of the general
results into one extended statement. Theorem~\ref{generalcase},
then, duplicates the essentials of our Abstract.

\begin{lemma}\label{L4.1}
Let $(X,\sT)$ be a crowded Tychonoff space. For $\emptyset\neq U\in\sT$
there is $V\in\sT$ such that $K:=\cl_{(X,\sT)}\,V$ satisfies
$V\subseteq K\subseteq U$ and
$\Delta(U)=\Delta(K)=|K|$.
\end{lemma}

Proof.
Choose $W\in\sT$ such that $W\subseteq U$ and $|W|=\Delta(U)$, and
choose $V\in\sT$ so that $V\neq\emptyset$ and $V\subseteq
K:=\cl_{(X,\sT)}\,V\subseteq W$.~$\square$

\begin{theorem}\label{generalcase}
Let $(X,\sT)$ be a crowded, maximally resolvable Tychonoff space such
that $S(X,\sT)\leq\Delta(X,\sT)=\kappa$. Then there are Tychonoff
expansions $\sU_i$ ($1\leq i\leq5$) of $\sT$, with
$\Delta(X,\sU_i)=\Delta(X,\sT)$ and $S(X,\sU_i)\leq\Delta(X,\sU_i)$,
such that $(X,\sU_i)$ is:

{\rm ($i=1$)}~$\omega$-resolvable but not
maximally resolvable;

{\rm ($i=2$)}~[if $\kappa'$ is regular, with
$S(X,\sT)\leq\kappa'\leq\kappa$] $\tau$-resolvable for all
$\tau<\kappa'$, but not
$\kappa'$-resolvable;

{\rm ($i=3$)}~maximally resolvable, but not extraresolvable;

{\rm ($i=4$)}~extraresolvable, but not maximally resolvable;

{\rm ($i=5$)}~maximally resolvable and extraresolvable, but not
strongly extraresolvable.
\end{theorem}

Proof. (Recall our frequently used convention that when $(X,\sT)$ is a
space and $Y\subseteq X$, the symbol $(Y,\sT)$ denotes  the set $Y$ with
the topology inherited from $(X,\sT)$.)

Using Lemma~\ref{L4.1} (with $U=X$), choose a regular-closed set $X'\subseteq X$
such that

$S(X',\sT)\leq S(X,\sT)\leq\Delta(X,\sT)=\Delta(X',\sT)=|X'|=\kappa$.

The definition of the topologies $\sU_i$ for $i=1,2,3$, and the
verification that they are as required, will be straightforward. We discuss
these first, leaving the cases $(i=4,5)$ for treatment later in the proof.

The space $(X',\sT)$ satisfies the hypotheses of Theorems \ref{maintheorem},
\ref{notkappa'resolv}, \ref{i=3}, so there are
Tychonoff expansions $\sU_i'$ ($i=1,2,3$) of $\sT$ on $X'$
satisfying their respective conclusions.
Let $\sU_i$ ($i=1,2,3$) be the topology on $X$ for which $(X',\sU_i')$
and $(X\backslash   X',\sT)$ are open-and-closed subspaces of
$(X,\sU_i)$.
It is easily seen that $(X,\sU_i)$ is a Tychonoff space. Further we have
$\sT\subseteq\sU_i$, since if $U\in\sT$ then $U\cap X'$ is open in
$(X',\sT)$, hence in $(X',\sU_i')$, hence in $(X',\sU_i)$, and
$U\cap(X\backslash X')$ is open in
$(X\backslash X',\sT)=(X\backslash X',\sU_i)$.

For $i=1,2,3$ we have, using
$\Delta(X',\sU_i')=\Delta(X',\sT)\leq\Delta(X\backslash X',\sT)$, that

$\Delta(X,\sU_i)=\min\{\Delta(X',\sU_i),\Delta(X\backslash X',\sU_i)\}
=\Delta(X',\sT)=\Delta(X,\sT)=\kappa$.

Further for $i=1,3$ we have, using $S(X',\sU_i')=S(X',\sT)$, that

$S(X,\sU_i)=S(X',\sU_i)+S(X\backslash X'\sU_i)
=S(X',\sT)+S(X\backslash X',\sT)=S(X,\sT)$,\\
\noindent while for $i=2$ we have

$S(X,\sU_2)=S(X',\sU_2)+S(X\backslash X',\sU_2)
=\kappa'+S(X\backslash X',\sU_2)=\kappa'+S(X\backslash X',\sT)=\kappa'$.

We verify the required (non-)
resolvability properties of the spaces
$(X,\sU_i)$ for $i=1,2,3$.

In each case, $(X,\sU_i)$ is the union of two disjoint open-and-closed
subspaces, namely $(X',\sU_i)$ and
$(X\backslash X',\sU_i)=(X\backslash X',\sT)$. When $i=1$,
these are both $\omega$-resolvable; when $i=2$,
both are $\tau$-resolvable for each $\tau<\kappa'$; when $i=3$, both are
$\kappa$-resolvable. Thus $(X,\sU_1)$ is $\omega$-resolvable;
$(X,\sU_2)$ is $\tau$-resolvable for all $\tau<\kappa'$; and $(X,\sU_3)$
is $\kappa$-resolvable (i.e., is maximally resolvable).

Since
$(X',\sU_1')=(X',\sU_1)$ is open in
$(X,\sU_1)$ and is not $\Delta(X,\sT)$-resolvable, surely $(X,\sU_1)$ is
not $\Delta(X,\sT)$-resolvable, i.e., is not $\Delta(X',\sU_1)$-resolvable.

The space $(X,\sU_2)$
is not $\kappa'$-resolvable, since its open subspace
$(X',\sU_2')=(X',\sU_2)$ is not
$\kappa'$-resolvable.

The space $(X,\sU_3)$ is not extraresolvable, since its open subspace
$(X',\sU_3')=(X',\sU_3)$ is not extraresolvable (and satisfies
$\Delta(X',\sU_3')=\Delta(X,\sU_3)$).

We turn to the cases
$(i=4,5)$.

Let $\sV\subseteq\sT$ be chosen maximal with respect to the properties

$\{\cl_{(X,\sT)}\,V:V\in \sV\}$ is pairwise disjoint, and

$|V|=|\cl_{(X,\sT)}\,V|=\Delta(V)$ for each $V\in\sV$.\\
\noindent We write $\sV=\{V_\beta:\beta<\alpha\}$ and
$X'_\beta:=\cl_{(X,\sT)}\,V_\beta$, the indexing chosen with $V_0$ and
$X_0'=X'$ as in the first part of this proof:
$|X_0'|=\Delta(X_0',\sT)=\Delta(X,\sT)$.

The space $(X_0',\sT)$ satisfies

$S(X_0',\sT)\leq
S(X,\sT)\leq\Delta(X,\sT)=\Delta(X_0',\sT)=|X_0|=\kappa$,\\
\noindent so by
Theorems~\ref{i=4}
and \ref{i=5} there are Tychonoff refinements
$\sU_{0,4}'$ and $\sU_{0,5}'$ of $(X_0',\sT)$, with

$S(X_0',\sU_{0,4}')=S(X_0',\sU_{0,5}')=S(X_0',\sT)$ and

$\Delta(X_0',\sU_{0,4}')=\Delta(X_0',\sU_{0,5}')=\Delta(X_0',\sT)=\kappa$,\\
\noindent
such that

$(X_0',\sU_{0,4}')$ is extraresolvable, but not maximally
resolvable; and

$(X_0',\sU_{0,5}')$ is maximally resolvable
and extraresolvable, but not strongly
extraresolvable.

For $0<\beta<\alpha$ the spaces $(X_\beta',\sT)$ satisfy

$S(X_\beta',\sT)\leq
S(X,\sT)\leq\kappa=\Delta(X,\sT)\leq\Delta(X_\beta',\sT)=|X_\beta'|$.\\
\noindent By Theorem~\ref{usingS_er}, taking
$\tau=\kappa_\beta:=|X_\beta'|$ there,
there are for $0<\beta<\alpha$
Tychonoff expansions $\sU_{\beta}'$ of
$(X_\beta',\sT)$ such that

$S(X_\beta',\sU_\beta')=S(X_\beta',\sT)$ and
$\Delta(X_\beta',\sU_\beta')=\Delta(X_\beta',\sT)$,\\
\noindent
and $(X_\beta',\sU_{\beta}')$ is
$\kappa_\beta$-resolvable and $2^{\kappa_\beta}$-extraresolvable. 
Then since $\kappa\leq\kappa_\beta$, the space
$(X_\beta',\sU_{\beta}')$ is
$\kappa$-resolvable and $2^{\kappa}$-extraresolvable.

Now for $(i=4,5)$ we define $\sU_i$ to be the smallest topology on $X$ such
that

(1) $\sT\subseteq\sU_i$,

(2) $(X'_0,\sU'_{0,i})$ is open-and-closed in
$(X,\sU_i)$, and

(3) each space $(X_\beta',\sU_{\beta}')$ (with $0<\beta<\alpha$)
is open-and-closed in
$(X,\sU_i)$.

To see that $(X,\sU_i)$ is a Tychonoff space, it is enough to note that
if $x\in\bigcup_{\beta<\alpha}\,X_\beta'$,
say $x\in X_{\overline{\beta}}'$, then
$X_{\overline{\beta}}'$ is an open Tychonoff neighborhood of
$x$ in $(X,\sU_i)$; while
if $x\notin\bigcup_{\beta<\alpha}\,X_\beta'$, then the $\sT$-open
neighborhoods of $x$ remain basic at $x$ in $(X,\sU_i)$ (so if $x\in
U\in\sU_i$ then there is a $\sU_i$-continuous (even, $\sT$-continuous)
real-valued function $f$ on $X$ such that $f(x)=0$ and $f=1$ on
$X\backslash U$).

For $\beta<\alpha$ we have

$\Delta(X_\beta',\sU_i)=\Delta(X_\beta',\sU_\beta')=\Delta(X_\beta',\sT)\geq\Delta(X_0',\sT)$,\\
\noindent so
$\Delta(X,\sU_i)=\min_{\beta<\alpha}\,\Delta(X_\beta',\sU_i)=\Delta(X_0',\sT)=\Delta(X,\sT)=\kappa$.

We verify for ($i=4,5$) that $S(X,\sU_i)\leq\Delta(X,\sU_i)$.
For a cellular family $\sW\subseteq\sU_i$ and $\beta<\alpha$ let
$\sW(\beta):=\{W\cap X_\beta':W\in\sW,W\cap X_\beta'\neq\emptyset\}$.
Then $\sW(\beta)$ is a cellular family
by Lemma~\ref{L4.1}. The set $\bigcup_{\beta<\alpha}\,X_\beta'$ is dense
in $(X,\sT)$, so $\sW=\bigcup_{\beta<\alpha}\,\sW(\beta)$, so
$|\sW|\leq\Sigma_{\beta<\alpha}\,|\sW(\beta)|$ with each

$|\sW(\beta)|<S(X_\beta',\sU_\beta')=S(X_\beta',\sT)\leq S(X,\sT)$.\\
\noindent Since $\alpha<S(X,\sT)$ and $S(X,\sT)$ is
regular, we have $|\sW|<S(X,\sT)$. It follows that

$S(X,\sU_i)\leq S(X,\sT)\leq\Delta(X,\sT)=\Delta(X,\sU_i)$.

It remains to verify that the spaces $(X,\sU_4)$ and $(X,\sU_5)$ have
the required (non-)
resolvability properties.

Each space $(X_\beta',\sU_4)$ is open in $(X,\sU_4)$, with
$\bigcup_{\beta<\alpha}\,X_\beta'$ dense in $(X,\sU_4)$. Each
space $(X_\beta',\sU_4)$ is extraresolvable (by Theorem~\ref{i=4}(b) for
$\beta=0$, by Theorem~\ref{usingS_er} for $0<\beta<\alpha$),
so for each $\beta<\alpha$ there is a family
$\sE_\beta=\{E_\beta(\eta):\eta<\kappa^+\}$ of dense subsets of
$(X_\beta',\sU_4)$ such that $E_\beta(\eta)\cap E_\beta(\eta')$ is
nowhere dense in $(X_\beta',\sU_4)$ whenever $\eta<\eta'<\kappa^+$. Then
with $E(\eta):=\bigcup_{\beta<\alpha}\,E_\beta(\eta)$, the family
$\{E(\eta):\eta<\kappa^+\}$ witnesses the extraresolvability of
$(X,\sU_4)$. The space $(X,\sU_4)$
is not maximally resolvable (i.e., is not $\kappa$-resolvable), however,
since its open subspace $(X_0',\sU_4)=(X_0',\sU_{4,0})$ is not
$\kappa$-resolvable.

The space $\bigcup_{\beta<\alpha}\,X_\beta'$ is open and dense in
$(X,\sU_5)$, with each $(X_\beta',\sU_5)$ open and $\kappa$-resolvable
and $2^\kappa$-extraresolvable, so $(X,\sU_5)$ is $\kappa$-resolvable
(i.e., maximally resolvable) and extraresolvable.
Each set $K\in[X_0']^{<\kappa}$ is closed and discrete in
$(X_0',\sU_{0,5})=(X_0',\sU_5)$, so
$\nwd(X_0',\sU_5)=\kappa$. Thus any family of sets dense in
$(X,\sU_5)$ witnessing the strong extraresolvability
of $(X,\sU_5)$ would trace on
$(X_0',\sU_5)$ to a family witnessing strong extraresolvability
there.~$\square$

\section{Some Questions}

Both our result cited from \cite{comfhu07} in
Remark~\ref{aboutkappa'-res}(b)
(where $S(X)>|X|$) and its sequel in Theorem~\ref{notkappa'resolv}(b)
(where $S(X)\leq|X|$) show that in some cases $\omega$-resolvability
suffices to guarantee $\tau$-resolvability for many larger $\tau$. Our
methods appear insufficiently delicate, however, to respond to the following
question.

\begin{question}\label{Q1}
{\rm
Let $(X,\sT)$ be an $\omega$-resolvable Tychonoff space such that
$S(X,\sT)\leq\Delta(X,\sT)$. Must $(X,\sT)$ be $\tau$-resolvable for
every $\tau<S(X,\sT)$?
}
\end{question}

\begin{question}\label{Q2}
{\rm
Let $X=(X,\sT)$ be a dense, $\omega$-resolvable subspace of the space
$(D(\kappa))^{2^\kappa}$ such that $|X|=\Delta(X)=\kappa$. [Then
$S(X)=\kappa^+$, and $X$ is
$\kappa$-resolvable, i.e., maximally resolvable, according to our result
\cite{comfhu07}(4.2).] Does $X$ admit a Tychonoff refinement $\sU$
(necessarily with $S(X,\sU)=\kappa^+$) such
that $\Delta(X,\sU)=\Delta(X,\sT)$, and $(X,\sU)$
is $\omega$-resolvable but not maximally resolvable? Always? Sometimes?
Never?
}
\end{question}

\begin{remarks}\label{aboutquestions}
{\rm
(a) Theorem \ref{maintheorem} sheds no light on Question~\ref{Q2}, since
the hypothesis $S(X,\sT)\leq\Delta(X,\sT)$ is lacking.

(b) The expansion $\sU$ of $\sT$ requested in Question~\ref{Q2}, if it
exists, cannot be of the kind constructed in this paper. More
specifically: There can
be no family $\sW\subseteq\sP(X)$ such that (i)~$|U\cap W|=\kappa$ for each
$W\in\sW$ and $\emptyset\neq U\in\sT$, (ii)~$\sU$ is the smallest
topology on $X$ containing $\sT$ and $\sW$, and (c)~each $W\in\sW$ is
$\sU$-clopen. For according to the argument outlined in
Discussion~\ref{disc1}, a space $(X,\sU)$ arising in that way will embed as a
dense subspace of $(D(\kappa))^I$ (with $|I|=w(X,\sU)$), hence if
$\omega$-resolvable is necessarily $\kappa$-resolvable.

(c) Many additional questions relating to (ir)resolvability, together
with extensive bibliographic citations, are recorded in the ``Problems"
article of Pavlov~\cite{pavlov07}.
}
\end{remarks}

\begin{remark}\label{nextpaper}
{\rm
The reader will have no difficulty using the methods of this paper
to establish the following result:

(*) {\it Let $(X,\sT)$ be a crowded,
maximally resolvable Tychonoff space with
$S(X,\sT)$ $\leq$ $\Delta(X, \sT)$ $=\kappa$. Then for
(fixed) $n<\omega$ there is
a Tychonoff expansion $\sU$ of $\sT$ such that $(X,\sU)$ is
$n$-resolvable but not $(n+1)$-resolvable.}

(Indeed, reducing as in Theorem~\ref{generalcase} to the case
$|X|=\Delta(X,\sT)$, it is enough to begin with a dense partition
$\{D^k_\eta:k<n,\eta<\kappa\}$ of $(X,\sT)$,
a strong small-set-separating
$\kappa$-independent family $\sI=\{\sI_t:t\in X\times2^\kappa\}$ of
$\kappa$ with each $\sI_t=\{I^0_t,I^1_t\}$, and with the family
$\sK=\{K_\xi:\xi<2^\kappa\}$ defined as in Theorem \ref{heredirr}.
Then the relation $X=\bigcup_{k<n}\,D^k$ with
$D^k:=\bigcup_{\eta<\kappa}\,D^k_\eta$ expresses $(X,\sU)$ with
$\sU:=\sT_{\mathcal{KID}}$ as the union of $n$-many disjoint dense sets,
each hereditarily irresolvable by Lemma~\ref{heredirr}(c). A space $(X,\sU)$
with such a partition cannot be
$(n+1)$-resolvable~\cite{illanes}, \cite{comfhu10?}.

We omit the details here of a proof of statement (*)
because a stronger Theorem is available, as
follows. 

(**) {\it For every $0<n<\omega$, every $n$-resolvable Tychonoff space
$(X,\sT)$ admits a Tychonoff expansion $\sU$ such that $(X,\sU)$ is
$n$-resolvable but not $(n+1)$-resolvable.}

We will prove (**)
in a manuscript now in preparation~\cite{comfhu08abs},
\cite{comfhu10?}. We remark {\it en
passant} that {\it ad hoc} constructions of Tychonoff spaces which for fixed
$n<\omega$ are $n$-resolvable but not $(n+1)$-resolvable have been
available for some time~\cite{vd93b}; see also
\cite{elkin69c}, \cite{fengmasa99}, \cite{eck97} and \cite{feng00}
for other examples, not all Tychonoff.

}
\end{remark}

\end{document}